\documentclass[11pt, reqno, psamsfonts]{amsart}
\pdfoutput=1

\usepackage{amssymb}
\usepackage{amsthm}
\usepackage{amsmath}
\usepackage{latexsym}
\usepackage[T1]{fontenc}
\usepackage[utf8]{inputenc}
\usepackage[russian, french, english]{babel}
\usepackage{graphicx}
\usepackage{wrapfig}
\usepackage[justification=centering, labelfont=bf]{caption}
\usepackage{mathtools}
\usepackage{tikz}
\usepackage{amsbsy}
\usepackage[inline]{enumitem}
\usepackage{mathrsfs}
\usetikzlibrary{shapes,snakes}
\usetikzlibrary{arrows.meta}
\usetikzlibrary{decorations.pathmorphing}
\usetikzlibrary{patterns}
\usepackage{float}
\usepackage{array}
\usepackage{multicol}
\usepackage{stmaryrd}
\usepackage{cancel} %to strike through text
\usepackage{lmodern}
\usepackage{mathabx}
\usepackage{titlesec}
\usepackage[spacing=true,kerning=true,babel=true,tracking=true]{microtype}
\usepackage[shortcuts]{extdash}
\usepackage[foot]{amsaddr}
\usepackage[left=1in,right=1in,top=1in,bottom=1in,bindingoffset=0cm]{geometry}
\usepackage{centernot}
\usepackage{mdframed}
\usepackage[hidelinks]{hyperref}
\usepackage{xspace}
\usepackage[backend=biber, style=alphabetic, sorting=nyt, maxnames=100,backref=true]{biblatex}
\title{\sffamily Equivariant maps to subshifts whose points have small stabilizers}
\date{}

\author{Anton~Bernshteyn}
\address{\normalfont School of Mathematics, Georgia Institute of Technology, Atlanta, GA, USA}
\email{bahtoh@gatech.edu}

\thanks{This research is partially supported by the NSF grant DMS-2045412.}

\newtheoremstyle{bfnote}%
{}{}%
{\slshape}{}%
{\bfseries}{\bfseries.}%
{ }%
{\thmname{#1}\thmnumber{ #2}\thmnote{ \ep{\normalfont{}#3}}}

\newtheoremstyle{claim}%
{}{}%
{\slshape}{}%
{\bfseries}{.}%
{ }%
{\thmname{#1}\thmnumber{ #2}\thmnote{ \ep{\normalfont{}#3}}}

\theoremstyle{bfnote}
\newtheorem{theo}{Theorem}[section]
\newtheorem*{theo*}{Theorem}
\newtheorem{prop}[theo]{Proposition}
\newtheorem{lemma}[theo]{Lemma}
\newtheorem{claim}[theo]{Claim}
\newtheorem*{claim*}{Claim}
\newtheorem{corl}[theo]{Corollary}

\newtheorem*{corl*}{Corollary}

\theoremstyle{definition}

\newtheorem*{defn*}{Definition}

\newtheorem{remk}[theo]{Remark}

\newtheorem*{exmp*}{Example}

\theoremstyle{remark}
\newtheorem*{ques*}{Question}
\newtheorem*{remk*}{Remark}

\theoremstyle{claim}
\newcounter{ForClaims}[section]
\newtheorem{subclaim}{Subclaim}[ForClaims]
\newtheorem{smallclaim}{Claim}[ForClaims]

\newcommand*{\myproofname}{Proof}
\newenvironment{claimproof}[1][\myproofname]{\begin{proof}[#1]}{\end{proof}}

\newcommand{\0}{\varnothing}
\newcommand{\set}[1]{\{#1\}}
\newcommand{\N}{{\mathbb{N}}}
\newcommand{\Z}{\mathbb{Z}}

\renewcommand{\P}{\mathbb{P}}

\renewcommand{\epsilon}{\varepsilon}

\renewcommand{\phi}{\varphi}
\renewcommand{\theta}{\vartheta}
\renewcommand{\leq}{\leqslant}
\renewcommand{\geq}{\geqslant}
\newcommand{\defeq}{\coloneqq}

\newcommand{\bemph}[1]{{\normalfont#1}} % define how emphasised brackets should look
\newcommand{\ep}[1]{\bemph{(}#1\bemph{)}} % parentheses

\newcommand{\emphd}[1]{{\fontseries{b}\selectfont\textsf{#1}}}
\newcommand{\acts}{\mathrel{\reflectbox{$\righttoleftarrow$}}}
\renewcommand{\G}{\Gamma}
\newcommand{\pto}{\dashrightarrow}
\newcommand{\dom}{\mathrm{dom}}
\newcommand{\Free}{\mathsf{Free}}

\newcommand{\rest}[2]{{{#1}\vert_{#2}}}
\newcommand{\Stab}{\mathrm{Stab}}
\newcommand{\Nbhd}{N}
\newcommand{\B}{\mathcal{B}}
\newcommand{\Col}{\mathsf{Col}}
\newcommand{\Ext}{\mathsf{Ext}}
\newcommand{\pr}{\mathsf{p}}
\newcommand{\de}{\mathsf{d}}
\newcommand{\vdeg}{\mathsf{vdeg}}
\newcommand{\ord}{\mathsf{ord}}
\newcommand{\pfun}[2]{[#1 \pto #2]}
\newcommand{\Sh}{\mathcal{S}}

\numberwithin{equation}{section}

\newenvironment{scproof}[1][]{\begin{proof}[\textsc{\upshape{Proof}}#1]}{\end{proof}}

\titleformat{\section}[block]{\large\bfseries\sffamily}{\thesection.}{1ex}{}
\titleformat{\subsection}[block]{\bfseries\sffamily}{\thesubsection.}{1ex}{}
\titleformat{\subsubsection}[block]{\itshape}{\bfseries\sffamily\upshape\thesubsubsection.}{1ex}{}

\titlespacing*{\section}{0pt}{*3}{*1}
\titlespacing*{\subsection}{0pt}{*3}{*1}
\titlespacing*{\subsubsection}{0pt}{*3}{*1}

\makeatletter
\newcommand{\neutralize}[1]{\expandafter\let\csname c@#1\endcsname\count@}
\makeatother

\newenvironment{propcopy}[1]
{\renewcommand{\thetheo}{\ref{#1}}%
	\neutralize{theo}\phantomsection
	\begin{prop}}
	{\end{prop}}
	
\newenvironment{lemmacopy}[1]
{\renewcommand{\thetheo}{\ref{#1}}%
	\neutralize{theo}\phantomsection
	\begin{lemma}}
	{\end{lemma}}

\renewbibmacro{in:}{}

\renewbibmacro*{volume+number+eid}{%
	\printfield{volume}%
	\setunit*{\addnbspace}% NEW (optional); there's also \addnbthinspace
	\printfield{number}%
	\setunit{\addcomma\space}%
	\printfield{eid}}
\DeclareFieldFormat[article]{volume}{\textbf{#1}\space}
\DeclareFieldFormat[article]{number}{\mkbibparens{#1}}

\DeclareFieldFormat{journaltitle}{#1,}
\DeclareFieldFormat[thesis]{title}{\mkbibemph{#1}\addperiod}
\DeclareFieldFormat[article, unpublished, thesis]{title}{\mkbibemph{#1},}
\DeclareFieldFormat[book]{title}{\mkbibemph{#1}\addperiod}
\DeclareFieldFormat[unpublished]{howpublished}{#1, }

\DeclareFieldFormat{pages}{#1}

\DeclareFieldFormat[article]{series}{Ser.~#1\addcomma}

\setlist{topsep=3pt,itemsep=3pt}

\usepackage{stackengine,amsmath}
\stackMath
\DeclareMathOperator*\dcap{{\stackinset{r}{-1.02ex}{c}{-1.9pt}{\downarrow}
  {\bigcap}\mkern2mu}}
\DeclareMathOperator*\acup{{\stackinset{r}{-1.02ex}{c}{1.9pt}{\uparrow}
  {\bigcup}\mkern2mu}}
\begin{document}

    \vspace*{-10pt}

    \maketitle
    
    %\vspace*{10pt}
    
    \begin{abstract}
        Let $\G$ be a countably infinite group. Given $k \in \N$, we use $\Free(k^\G)$ to denote the free part of the Bernoulli shift action $\G \acts k^\G$. Seward and Tucker-Drob showed that there exists a free subshift $\Sh \subseteq \Free(2^\G)$ such that every free Borel action of $\G$ on a Polish space admits a Borel $\G$-equivariant map to $\Sh$. Here we generalize this result as follows. Let $\Sh$ be a subshift of finite type \ep{for example, $\Sh$ could be the set of all proper colorings of the Cayley graph of $\G$ with some finite number of colors}. Suppose that $\pi \colon \Free(k^\G) \to \Sh$ is a continuous $\G$-equivariant map and let $\Stab(\pi)$ be the set of all group elements that fix every point in the image of $\pi$. Unless $\pi$ is constant, $\Stab(\pi)$ is a finite normal subgroup of $\G$. We prove that there exists a subshift $\Sh' \subseteq \Sh$ such that:
        \begin{itemize}
	        \item the stabilizer of every point in $\Sh'$ is $\Stab(\pi)$, and
	        \item every free Borel action of $\G$ on a Polish space admits a Borel $\G$-equivariant map to $\Sh'$.
	    \end{itemize}
	    %Taking $X = 2^\G$ and letting $\pi$ be the inclusion $\Free(2^\G) \to 2^\G$, we recover the Seward--Tucker-Drob theorem. Our proof uses probabilistic tools developed earlier by the present author.
	    %This generalizes a result of Seward and Tucker-Drob. %, who showed that the above statement holds when $X = 2^\G$ and $\pi$ is the inclusion map $\Free(2^\G) \to 2^\G$. %We then apply our result to proper colorings of Schreier graphs.
	    In particular, %if $\Stab(\pi)$ is trivial, i.e.,
	    if the shift action of $\G$ on the image of $\pi$ is faithful (i.e., if $\Stab(\pi)$ is trivial), then the subshift $\Sh'$ is free. As an application of this general result, we deduce that if $F$ is a %nonempty
	    finite symmetric subset of $\G \setminus \set{\mathbf{1}}$ of size $|F| = d \geq 1$ %not containing the identity
	    and $\Col(F, d + 1) \subseteq (d+1)^\G$ is the set of all proper $(d+1)$-colorings of the Cayley graph of $\G$ corresponding to $F$, then there is a free subshift $\Sh \subseteq \Col(F, d+1)$ such that every free Borel action of $\G$ on a Polish space admits a Borel $\G$-equivariant map to $\Sh$.
    \end{abstract}
	
	\section{Introduction}
	
	Throughout this paper, $\G$ is a countably infinite group with identity element $\mathbf{1}$. Given an integer $k \in \N$, we identify $k$ with the $k$-element set $\set{0, 1, \ldots, k-1}$ and equip it with the discrete topology. The \emphd{\ep{Bernoulli} shift} is the action of $\G$ on the product space $k^\G$ given by the formula
	\[
	    (\gamma \cdot x)(\delta) \,\defeq\, x(\delta \gamma) \quad \text{for all } x \colon \G \to k \text{ and } \gamma, \, \delta \in \G.
	\]
	The product topology on $k^\G$ is compact, metrizable, and zero-dimensional \ep{that is, it has a base consisting of clopen sets}, and the shift action $\G \acts k^\G$ is continuous. The \emphd{free part} of $k^\G$ is the set
	\[\Free(k^\G) \,\defeq\, \set{x \in k^\G \, : \, \Stab(x) = \set{\mathbf{1}}},\]
	where $\Stab(x)$ denotes the stabilizer of $x$. In other words, $\Free(k^\G)$ is the largest subspace of $k^\G$ on which the shift action is free. %Note that $\Free(k^\G)$ is a $G_\delta$ subset of $k^\G$, and hence it is a Polish space \cite[Theorem 3.11]{KechrisDST}.
	If $k \geq 2$, then $\Free(k^\G)$ is dense in $k^\G$ \ep{see, e.g., \cite[Lemma 2.3]{Bernoulli}}.
	A \emphd{subshift} is a closed shift-invariant subset $\Sh \subseteq k^\G$. A subshift $\Sh \subseteq k^\G$ is \emphd{free} if $\Sh \subseteq \Free(k^\G)$. %, i.e., if $\G$ acts freely on $X$.
	An important class of subshifts are \emphd{subshifts of finite type}, or \emphd{SFTs} for short, which are of the form %\[X \,=\, \bigcap \set{\gamma \cdot C \,:\, \gamma \in \G}\]
	\[%\begin{equation}\label{eq:clop}
	    \Sh \,=\, \bigcap_{\gamma \in \G} (\gamma \cdot C)
	\]%\end{equation}
	for some clopen $C \subseteq k^\G$. Equivalently, $\Sh$ is an SFT if and only if there exist a finite set $W \subset \G$, called a \emphd{window}, and a collection $\Phi \subseteq k^W$ of maps $\phi \colon W \to k$ such that
	\[
	    \Sh \,=\, \set{x \in k^\G \,:\, \rest{(\gamma \cdot x)}{W} \in \Phi \text{ for all } \gamma \in \G},
	\]
	where $\rest{(\gamma \cdot x)}{W}$ denotes the restriction of the function $\gamma \cdot x$ to $W$.
	
	The starting point of our investigation is the following result of Seward and Tucker-Drob:
	
	\begin{theo}[{Seward--Tucker-Drob \cite{STD}}]\label{theo:STD}
	    There exists a free subshift $\Sh \subset 2^\G$ such that every free Borel action $\G \acts X$ on a Polish space admits a Borel $\G$-equivariant map $X \to \Sh$. 
	\end{theo}
	
	Even the {existence} of a nonempty free subshift is far from obvious. For several special classes of groups, it was established by Dranishnikov and Schroeder \cite{DS} and Glasner and Uspenskij \cite{GU}, while the general case was settled by Gao, Jackson, and Seward \cite{GJS1, GJS2}. An alternative probabilistic construction was subsequently discovered by Aubrun, Barbieri, and Thomass\'e \cite{ABT}. Seward and Tucker-Drob's proof of Theorem~\ref{theo:STD} further develops the methods of \cite{GJS1, GJS2}. A simple probabilistic proof of Theorem~\ref{theo:STD} is given in \cite{Ber_cont}.
	
	Here we study $\G$-equivariant maps from free Borel actions of $\G$ to SFTs. It follows from Theorem~\ref{theo:STD} that if $\Sh$ is an SFT and there is a Borel $\G$-equivariant map $\pi \colon \Free(2^\G) \to \Sh$, then every free Borel action $\G \acts X$ on a Polish space admits a Borel $\G$-equivariant map to $\Sh$. Indeed, by Theorem~\ref{theo:STD}, there is a Borel $\G$-equivariant map $X \to \Free(2^\G)$; composing it with $\pi$ yields a desired map $X \to \Sh$. In \cite{Ber_cont} it was shown that an analogous statement holds for {continuous} $\G$-equivariant maps % The following fact is proved in \cite{Ber_cont}
	\ep{this result was also established independently by Seward; see \cite[\S7.1]{Toast}}:
	
	\begin{theo}[{\cite[Theorem~1.12]{Ber_cont}}]\label{theo:cont}
	    Let $\Sh$ be an SFT and suppose that there is a continuous $\G$-equivariant map $\Free(2^\G) \to \Sh$. Then every free continuous action $\G \acts X$ on a zero-dimensional Polish space admits a continuous $\G$-equivariant map $X \to \Sh$.
	\end{theo}
	
    Our aim is to strengthen Theorem~\ref{theo:STD} as follows: Given an SFT $\Sh$ and a continuous $\G$-equivariant map $\Free(k^\G) \to \Sh$, we wish to find a subshift $\Sh' \subseteq \Sh$ such that \begin{enumerate*}[label=\ep{\emph{\alph*}}]
        \item every free Borel action of $\G$ on a Polish space admits a Borel $\G$-equivariant map to $\Sh'$, and
        \item\label{item:almostfree} $\Sh'$ is as close to being free as possible.
    \end{enumerate*}
    To state item \ref{item:almostfree} precisely, we introduce the following notation. Let $\G \acts X$, $\G \acts Y$ be actions of $\G$ and let $\pi \colon X \to Y$ be a $\G$-equivariant map. Define %for a $\G$-equivariant map $\pi \colon \Free(k^\G) \to \ell^\G$, let
	\[
	    \Stab(\pi) \,\defeq\, \bigcap \set{\Stab(\pi(x)) \,:\, x \in X}.
	\]
	In other words, $\Stab(\pi)$ comprises the group elements that fix every point in the image of $\pi$. In particular, $\Stab(\pi)$ is trivial if and only if the action of $\G$ on $\pi(X)$ is \emphd{faithful}, meaning that the only group element $\gamma \in \G$ with $\gamma \cdot y = y$ for all $y \in \pi(X)$ is $\gamma = \mathbf{1}$. %Our goal is to find $X$ as in Theorem~\ref{theo:STD} such that $\Stab(\pi(x)) = \Stab(\pi)$ for all $x \in X$.
	%Let us record the following observation:
	
	\begin{prop}\label{prop:Stab}
	    Let $\pi \colon \Free(k^\G) \to \ell^\G$ be a continuous $\G$-equivariant map. Then either $\pi$ is constant and $\Stab(\pi) = \G$, or else, $\Stab(\pi)$ is a finite normal subgroup of $\G$.
	\end{prop}
	
	The \ep{easy} proof of Proposition~\ref{prop:Stab} is given in \S\ref{subsec:stab}. Note that if $\G$ has no nontrivial finite normal subgroups \ep{for example, if every non-identity element of $\G$ has infinite order or if $\G$ is ICC} and $\pi \colon \Free(k^\G) \to \ell^\G$ is a non-constant continuous $\G$-equivariant map, then we must have $\Stab(\pi) = \set{\mathbf{1}}$. In general, $\Stab(\pi)$ can be any finite normal subgroup of $\G$. To see this, let $H \lhd \G$ be a finite normal subgroup and let %Indeed, for a finite normal subgroup $H \lhd \G$, let
	$\mathcal{X}_H \subseteq 2^\G$ be the set of all functions $x \colon \G \to 2$ that are constant on the cosets of $H$. \ep{If $H = \set{\mathbf{1}}$, then $\mathcal{X}_H = 2^\G$.} Then $\mathcal{X}_H$ is an SFT and we have the following:  %here is also a converse to Proposition~\ref{prop:Stab}:
	
	\begin{prop}\label{prop:Stab_conv}
	    For every finite normal subgroup $H \lhd \G$ and every $k \geq 2$, there is a continuous $\G$-equivariant map $\pi \colon \Free(k^\G) \to \mathcal{X}_H$ such that $\Stab(\pi) = H$.
	\end{prop}
	
	The proof of Proposition~\ref{prop:Stab_conv} is also given in \S\ref{subsec:stab}. We are now ready to state our main result.
	
	\begin{theo}[Subshifts with small stabilizers]\label{theo:main}
	    Let $\Sh$ be an SFT and let $\pi \colon \Free(k^\G) \to \Sh$ be a continuous $\G$-equivariant map for some $k \geq 2$. Then there exists a subshift $\Sh' \subseteq \Sh$ such that:
	    \begin{itemize}
	        \item the stabilizer of every point in $\Sh'$ is $\Stab(\pi)$, and
	        \item every free Borel action $\G \acts X$ on a Polish space admits a Borel $\G$-equivariant map $X \to \Sh'$.
	    \end{itemize}
	\end{theo}
	
	As the case $\Stab(\pi) = \set{\mathbf{1}}$ is of particular interest, we record it as a corollary for ease of reference:
	
	\begin{corl}[Free subshifts from faithful actions]\label{corl:mainfree}
	    Let $\Sh$ be an SFT and let $\pi \colon \Free(k^\G) \to \Sh$ be a continuous $\G$-equivariant map for some $k \geq 2$. If the shift action of $\G$ on the image of $\pi$ is faithful, then there exists a free subshift $\Sh' \subseteq \Sh$ such that every free Borel action $\G \acts X$ on a Polish space admits a Borel $\G$-equivariant map $X \to \Sh'$ \ep{in particular, $\Sh' \neq \0$}.
	\end{corl}
	
	Theorem~\ref{theo:STD} follows by applying Corollary~\ref{corl:mainfree} %Theorem~\ref{theo:main}
	to $\Sh = 2^\G$ and the inclusion map $\iota \colon \Free(2^\G) \to 2^\G$, since the action $\G \acts \Free(2^\G)$ is, by definition, free and hence faithful. %$\Stab(\iota) = \set{\mathbf{1}}$.
	More generally, if $H \lhd \G$ is any finite normal subgroup, then, by Proposition~\ref{prop:Stab_conv}, there is a continuous $\G$-equivariant map $\pi \colon \Free(2^\G) \to \mathcal{X}_H$ with $\Stab(\pi) = H$. Theorem~\ref{theo:main} then yields a subshift $\mathcal{X}_H' \subseteq \mathcal{X}_H$ such that \begin{enumerate*}[label=\ep{\emph{\alph*}}]
        \item every free Borel action of $\G$ on a Polish space admits a Borel $\G$-equivariant map to $\mathcal{X}_H'$, and
        \item the stabilizer of every point in $\mathcal{X}_H'$ is $H$---which is as small as it can be since, by construction, $H$ fixes every point in $\mathcal{X}_H$.
    \end{enumerate*}
    %Note that every point in $\mathcal{X}_H$ is fixed by $H$, so the second
	
	In general, finding free subshifts with desired properties can be a subtle and challenging problem. It is therefore notable that Corollary~\ref{corl:mainfree} provides a general sufficient condition for the existence of a nonempty free subshift in a given SFT. %(moreover, not only is the subshift it yields nonempty, but it also has the property that every free Borel action of $\G$ on a Polish space admits a Borel $\G$-equivariant map to it).
	As an illustration, we describe some combinatorial consequences of Corollary~\ref{corl:mainfree}. %Theorem~\ref{theo:main}.
	Let $F \subset \G$ be a finite set not containing $\mathbf{1}$. Assume that $F$ is \emphd{symmetric}, i.e., $F = F^{-1}$. The \emphd{Cayley graph} $G(\G, F)$ is the graph with vertex set $\G$ in which two group elements $\gamma$, $\delta$ are adjacent if and only if $\delta = \sigma \gamma$ for some $\sigma \in F$. For $\ell \in \N$, a \emphd{proper $\ell$-coloring} of $G(\G, F)$ is a function $x \colon \G \to \ell$ such that $x(\gamma) \neq x(\delta)$ whenever $\gamma$ and $\delta$ are adjacent in $G(\G, F)$. Let $\Col(F, \ell) \subseteq \ell^\G$ be the set of all proper $\ell$-colorings of $G(\G, F)$. Then $\Col(F, \ell)$ is an SFT. If $\ell \geq |F| + 1$, then it is easy to construct a continuous $\G$-equivariant map $\pi \colon \Free(3^\G) \to \Col(F, \ell)$ such that $\Stab(\pi) = \set{\mathbf{1}}$; see \S\ref{subsec:col} for details. Corollary~\ref{corl:mainfree} then yields the following conclusion:
	
	\begin{corl}[Free subshift of proper colorings]\label{corl:coloring}
	    Let $F \subset \G$ be a nonempty finite symmetric set with $\mathbf{1} \not \in F$. If $\ell \geq |F| + 1$, then there exists a free subshift $\Sh \subseteq \Col(F, \ell)$ such that every free Borel action $\G \acts X$ on a Polish space admits a Borel $\G$-equivariant map $X \to \Sh$.
	\end{corl}
	
	The bound on $\ell$ in Corollary~\ref{corl:coloring} is, in general, optimal. For instance, if $H$ is a finite subgroup of $\G$ and $F = H \setminus \set{\mathbf{1}}$, then $\Col(F, |F|)$ is empty, because for any $x \in \Col(F, |F|)$, the values $x(h)$, $h \in H$ have to be distinct, and there are $|H| = |F| + 1$ of them. A more interesting example is due to Marks \cite{Marks}, who showed that if $\G = \Z_2^{\ast n}$ is the free product of $n$ copies of $\Z_2$ and $F$ is the standard generating set for $\G$, then $|F| = n$ and there is no Borel $\G$-equivariant map $\Free(2^\G) \to \Col(F, n)$.
	
	On the other hand, for some choices of $\G$ and $F$, the bound on $\ell$ can be improved. For instance, suppose that $\G = \Z^n$ for $n \geq 2$ and $F$ is the standard symmetric generating set for $\Z^n$ of size $2n$. In this setting, Gao and Jackson \cite[Theorem 4.2]{GaoJack} constructed a continuous $\Z^n$-equivariant map $\pi \colon \Free(2^{\Z^n}) \to \Col(F,4)$. Since $\Z^n$ has no nontrivial finite subgroups, $\Stab(\pi)$ must be trivial, so Corollary~\ref{corl:mainfree} yields a free subshift $\Sh \subseteq \Col(F, 4)$ such that every free Borel action $\Z^n \acts X$ on a Polish space admits a Borel $\Z^n$-equivariant map $X \to \Sh$. Actually, Chandgotia and Unger \cite{CU} showed that in this case, the optimal value for $\ell$ is $3$. However, the reduction from $4$ to $3$ cannot be achieved using Theorem~\ref{theo:main},  %Chandgotia and Unger further showed every free Borel action $\Z^d \acts Y$ on a Polish space admits a Borel $\Z^d$-equivariant map $\sigma \colon Y \to \Col(F,3)$ such that $\Z^d$ acts freely on $\sigma(Y)$. This result cannot be deduced from Theorem~\ref{theo:main},
	since, according to a result of Gao, Jackson, Krohne, and Seward \cite[Theorem 4.3]{Abelian}, there is no continuous $\Z^n$-equivariant map $\Free(2^{\Z^n}) \to \Col(F, 3)$.
	
	Let us now say a few words about the proof of Theorem~\ref{theo:main}.
	Curiously, an important role in it is played by subshifts of the form $\Col(F, \ell)$. % play an important role in the proof of Theorem~\ref{theo:main}.
	The argument consists of two lemmas. %We split the proof into two stages, %Let us now say a few words about the proof of Theorem~\ref{theo:main}. The proof is split into two steps 
	%encapsulated in the following lemmas: %, proved in \S\S\ref{sec:cont+Stab} and \ref{sec:col+Stab}:
	
	\begin{lemma}\label{lemma:cont+Stab}
	       Let $\Sh$ be an SFT and let $\pi \colon \Free(k^\G) \to \Sh$ be a continuous $\G$-equivariant map for some $k \geq 2$. Then there exist a finite symmetric set $F \subset \G$ with $\mathbf{1} \not \in F$, an integer $\ell \geq |F| + 1$, and a continuous $\G$-equivariant map $\tilde{\pi} \colon \Col(F,\ell) \to \Sh$ such that $\Stab(\tilde{\pi}) = \Stab(\pi)$.
	\end{lemma}
	
	\begin{lemma}\label{lemma:col+Stab}
	       Let $F \subset \G$ be a finite symmetric set with $\mathbf{1} \not \in F$. Fix integers $\ell \geq |F| + 1$, $m \geq 1$ and let $\pi \colon \Col(F,\ell) \to m^\G$ be a continuous $\G$-equivariant map. Then there exists a subshift $\mathcal{Z} \subseteq  \Col(F,\ell)$ with the following properties:
	    \begin{itemize}
	        \item for all $z \in \mathcal{Z}$, $\Stab(\pi(z)) = \Stab(\pi)$, and
	        \item every free Borel action $\G \acts X$ on a Polish space admits a Borel $\G$-equivariant map $X \to \mathcal{Z}$.
	    \end{itemize}
	\end{lemma}
	
	%Lemmas~\ref{lemma:cont+Stab} and \ref{lemma:col+Stab} are proved in \S\S\ref{sec:cont+Stab} and \ref{sec:col+Stab} respectively.
	Lemma~\ref{lemma:cont+Stab} allows us to replace $\Free(k^\G)$ with $\Col(F,\ell)$. One of the key advantages of working with the space $\Col(F, \ell)$ is that it is compact, and its compactness plays a crucial role in the proof of Lemma~\ref{lemma:col+Stab}. The main tool we rely on to establish Lemmas~\ref{lemma:cont+Stab} and \ref{lemma:col+Stab} is the continuous version of the Lov\'asz Local Lemma \ep{LLL} developed by the author in \cite{Ber_cont}, which we describe in \S\ref{sec:LLL}. In \cite{Ber_cont}, the continuous LLL was applied to prove Theorem~\ref{theo:cont} and to give a new simple proof of Theorem~\ref{theo:STD}. The proofs of Lemmas~\ref{lemma:cont+Stab} and~\ref{lemma:col+Stab} follow a similar outline but are significantly more technically involved. They are presented in \S\ref{sec:proof}. % respectively. %Both of them are obtained by extending and carefully modifying the approach used to prove Theorems \ref{theo:STD} and~\ref{theo:cont} in \cite[\S4]{Ber_cont}.
	
	Once Lemmas~\ref{lemma:cont+Stab} and \ref{lemma:col+Stab} have been verified, they quickly yield Theorem~\ref{theo:main}:
	
	\begin{scproof}[ of Theorem~\ref{theo:main}]%Indeed,
	Let $\Sh$ be an SFT and let $\pi \colon \Free(k^\G) \to \Sh$ be a continuous $\G$-equivariant map for some $k \geq 2$. By Lemma~\ref{lemma:cont+Stab}, there exist a finite symmetric set $F \subset \G$ with $\mathbf{1} \not \in F$, an integer $\ell \geq |F| + 1$, and a continuous $\G$-equivariant map $\tilde{\pi} \colon \Col(F,\ell) \to \Sh$ such that $\Stab(\tilde{\pi}) = \Stab(\pi)$. Now, by Lemma~\ref{lemma:col+Stab}, there is a subshift $\mathcal{Z} \subseteq \Col(F,\ell)$ such that:
	\begin{itemize}
	        \item for all $z \in \mathcal{Z}$, $\Stab(\tilde{\pi}(z)) = \Stab(\tilde{\pi}) = \Stab(\pi)$; and
	        \item every free Borel action $\G \acts X$ on a Polish space admits a Borel $\G$-equivariant map $X \to \mathcal{Z}$.
	\end{itemize}
	Taking $\Sh' \defeq \tilde{\pi}(\mathcal{Z})$ finishes the proof.
	\end{scproof}
	
	%As mentioned earlier, in \cite{Ber_cont} the author presented a new proof of the \hyperref[theo:STD]{Seward--Tucker-Drob theorem} based on probabilistic ideas. Our proof of Theorem~\ref{theo:main} follows a similar outline but is somewhat more technically involved, as we have to take special care to control both $\Stab(x)$ and $\Stab(\pi(x))$ for each $x \in X$.
	
	%\subsubsection*{{{Acknowledgments}}}
	
	\emph{Acknowledgments}.---I am very grateful to Marcin Sabok and Anush Tserunyan for providing a stimulating and productive environment during the \ep{virtual} \emph{Logic and Applications} section of the 2020 CMS Winter Meeting. This project arose from a question asked by Clinton Conley during Spencer Unger's talk there, namely whether Corollary~\ref{corl:coloring} holds. I am also greateful to the anonymous referees for their helpful feedback.
	
	\section{Preliminaries}\label{sec:prelim}
	
	\subsection{Basics of continuous combinatorics}\label{subsec:cont_graph}
	
	Let $G$ be a graph. For a subset $S \subseteq V(G)$, $\Nbhd_G(S)$ denotes the \emphd{neighborhood} of $S$ in $G$, i.e., the set of all vertices that have a neighbor in $S$. For a vertex $x \in V(G)$, we write $\Nbhd_G(x) \defeq \Nbhd_G(\set{x})$. A graph $G$ is \emphd{locally finite} if $\Nbhd_G(x)$ is finite for every $x \in V(G)$. The \emphd{maximum degree} $\Delta(G)$ of a graph $G$ is defined by $\Delta(G) \defeq \sup_{x \in V(G)} |\Nbhd_G(x)|$. A set $I \subseteq V(G)$ is \emphd{independent} if $I \cap \Nbhd_G(I) = \0$, i.e., if no two vertices in $I$ are adjacent. For a subset $U \subseteq V(G)$, we use $G[U]$ to denote the \emphd{subgraph of $G$ induced by $U$}, i.e., the graph with vertex set $U$ whose adjacency relation is inherited from $G$, and we write $G - U \defeq G [V(G) \setminus U]$.
	
	We say that a graph $G$ is \emphd{continuous} if $V(G)$ is a zero-dimensional Polish space and for every clopen set $U \subseteq V(G)$, its neighborhood $\Nbhd_G(U)$ is also clopen. Note that if $G$ is a continuous graph and $U \subseteq V(G)$ is a clopen set of vertices, then the subgraph $G[U]$ of $G$ induced by $U$ is also continuous. The following basic facts are standard. In the Borel context, they were proved by Kechris, Solecki, and Todorcevic in their seminal paper \cite{KST}. %The same proofs give clopen, rahter than just Borel, constructions.  %, where their analogs are proved in the Borel context. %Since their proofs are very short, we include them here for completeness.
	
	\begin{lemma}[{\cite[Lemma 2.1]{Ber_cont}}]\label{lemma:countable}
		Every locally finite continuous graph $G$ admits a partition $V(G) = \bigsqcup_{n = 0}^\infty I_n$ into countably many clopen independent sets.
	\end{lemma}
	%\begin{scproof}
	%	Let $\set{U_n \,:\, n \in \N}$ be a countable base for the topology on $V(G)$ consisting of clopen sets. For each $n \in \N$, let $V_n \defeq U_n \setminus \Nbhd_G(U_n)$. By construction, each $V_n$ is independent and, since $G$ is continuous, clopen. Since $G$ is locally finite, each $x \in V(G)$ has an open neighborhood disjoint from $\Nbhd_G(x)$, and hence $\bigcup_{n = 0}^\infty V_n = V(G)$. It remains to make the sets disjoint by setting $I_n \defeq V_n \setminus (V_0 \cup \ldots \cup V_{n-1})$.
	%\end{scproof}

	\begin{lemma}\label{lemma:max}
		Let $G$ be a locally finite continuous graph and let $J \subseteq V(G)$ be a clopen independent set. Then there is a clopen maximal independent set $I \subseteq V(G)$ with $I \supseteq J$.
	\end{lemma}
	\begin{scproof}
		Let $V(G - J) = \bigsqcup_{n = 1}^\infty I_n$ be a partition into countably many clopen independent sets given by Lemma~\ref{lemma:countable}. Define a sequence of clopen subsets $I_n' \subseteq I_n$ recursively by setting $I_0' \defeq J$ and $I_{n+1}' \defeq I_{n+1} \setminus \Nbhd_G(I_0' \sqcup \ldots \sqcup I_n')$ for all $n \in \N$. We claim that the set $I \defeq \bigsqcup_{n=0}^\infty I_n'$ is as desired. Indeed, by construction, $I$ is a maximal independent set containing $J$. Since $G$ is continuous, the sets $I_n'$ are clopen, and hence $I$ is open. But the sets $I_n \setminus I_n'$ are also clopen, so $V(G) \setminus I = \bigsqcup_{n=0}^\infty (I_n \setminus I_n')$ is open as well, and hence $I$ is clopen, as claimed. 
	\end{scproof}
	
	Lemma~\ref{lemma:max} is often applied with $J = \0$, in which case it simply asserts the existence of a clopen maximal independent set $I \subseteq V(G)$.
	
	Recall that a \emphd{proper $\ell$-coloring} of a graph $G$ is a function $f \colon V(G) \to \ell$ such that $f(x) \neq f(y)$ whenever $x$ and $y$ are adjacent vertices of $G$.

	\begin{lemma}\label{lemma:coloring}
		Let $G$ be a continuous graph of finite maximum degree $\Delta$. If $\ell \geq \Delta + 1$ and $C \subseteq V(G)$ is a clopen set, then every continuous proper $\ell$-coloring $g \colon C \to \ell$ of $G[C]$ can be extended to a continuous proper $\ell$-coloring $f \colon V(G) \to \ell$ of $G$. %is a continuous proper $\ell$-coloring of $G[C]$, then  and let $J_0$, \ldots, $J_{\ell-1} \subseteq V(G)$ be disjoint clopen independent sets. Then $G$ admits a partition $V(G) = I_0 \sqcup \ldots \sqcup I_{\ell-1}$ into $\ell$ clopen independent sets such that $I_i \supseteq J_i$ for each $i$.
	\end{lemma}
	\begin{scproof}
		For $i < \ell$, let $J_i \defeq g^{-1}(i)$, so $C = J_0 \sqcup \ldots \sqcup J_{\ell - 1}$ is a partition into clopen independent sets. Apply Lemma~\ref{lemma:max} iteratively to obtain a sequence $I_0$, \ldots, $I_{\ell-1}$, where each $I_i$ is a clopen maximal independent set in the graph $G - I_0 - \cdots - I_{i-1} - J_{i+1} - \cdots - J_{\ell-1}$ and $I_i \supseteq J_i$. We claim that $V(G) = I_0 \sqcup \ldots \sqcup I_{\ell-1}$. Indeed, every vertex not in $I_0 \sqcup \ldots \sqcup I_{\ell-1}$ must have a neighbor in each of $I_0$, \ldots, $I_{\ell-1}$, which is impossible as the maximum degree of $G$ is at most $\ell-1$. Setting $f(x)$ to be the unique index $i$ such that $x \in I_i$ finishes the proof.
	\end{scproof}
	
	%Lemma~\ref{lemma:coloring} is often applied with $J_0 = \cdots = J_\Delta = \0$, in which case it just says that $V(G)$ can be partitioned into $\Delta + 1$ clopen independent sets.
	
	\subsection{Schreier graphs and coding maps}\label{subsec:code}
	
	Let $\G \acts X$ be a continuous action of $\G$ on a zero-dimensional Polish space $X$. Given a finite set $F \subset \G$, we say that the action $\G \acts X$ is \emphd{$F$-free} if $F \cap \Stab(x) \subseteq \set{\mathbf{1}}$ for all $x \in X$, i.e., if every non-identity element of $F$ acts freely on $X$. The \emphd{Schreier graph} $G(X,F)$ is the graph with vertex set $X$ in which two distinct vertices $x$, $y$ are adjacent if and only if $y = \sigma \cdot x$ for some $\sigma \in F \cup F^{-1}$. If the action $\G \acts X$ is $F$-free, then for a clopen set $U \subseteq X$, its neighborhood in $G(X,F)$ is $((F \cup F^{-1})\setminus \set{\mathbf{1}}) \cdot U$, which is also clopen. Therefore, in this case $G(X,F)$ is continuous.
	%No generality is lost by assuming that $F$ is symmetric and $\mathbf{1} \not \in F$ since, by definition, \[G(X,F) \,=\, G(X, (F \cup F^{-1})\setminus \set{\mathbf{1}}).\] 
	Note that the Cayley graph $G(\G, F)$ is a special case of this construction applied to the left multiplication action $\G \acts \G$ of $\G$ on itself \ep{viewed as a discrete space}.
	
	For an action $\G \acts X$, there is a natural one-to-one correspondence
	\[
	    \set{\text{functions $X \to k$}} \quad \longleftrightarrow \quad \set{\text{$\G$-equivariant maps $X \to k^\G$}}.
	\]
	Namely, each function $f \colon X \to k$ gives rise to the $\G$-equivariant \emphd{coding map} $\pi_f \colon X \to k^\G$ given by
	\[
	    \pi_f(x)(\gamma) \,\defeq\, f(\gamma \cdot x) \quad \text{for all } x \in X \text{ and } \gamma \in \G.
	\]
	Conversely, if $\pi \colon X \to k^\G$ is $\G$-equivariant, then $\pi = \pi_f$ for the function $f= (x \mapsto \pi(x)(\mathbf{1}))$. As a special case of this correspondence, observe that proper $\ell$-colorings of the Schreier graph $G(X, F)$ correspond exactly to $\G$-equivariant maps $X \to \Col(F, \ell)$.
	
	\subsection{Stabilizers and finite normal subgroups}\label{subsec:stab}
	
	Here we prove Propositions~\ref{prop:Stab} and \ref{prop:Stab_conv} from the introduction.
	
	\begin{propcopy}{prop:Stab}
	    Let $k$, $\ell \geq 1$ and let $\pi \colon \Free(k^\G) \to \ell^\G$ be a continuous $\G$-equivariant map. Then either $\pi$ is constant and $\Stab(\pi) = \G$, or else, $\Stab(\pi)$ is a finite normal subgroup of $\G$.
	\end{propcopy}
	\begin{scproof}
	    Assume that $\pi$ is not constant. It is clear that $\Stab(\pi)$ is a normal subgroup of $\G$, so we only need to show that $\Stab(\pi)$ is finite. To this end, let
	    \[f \colon \Free(k^\G) \to \ell \colon x \mapsto \pi(x)(\mathbf{1}).\] \ep{In other words, $f$ is such that $\pi = \pi_f$.} Since $\pi$ is not constant, there are $x$, $y \in \Free(k^\G)$ such that $f(x) \neq f(y)$. Since $f$ is continuous, it is constant on some open neighborhoods of $x$ and $y$. This means that there is a finite set $W \subset \G$ such that for all $z \in \Free(k^\G)$, $f(z) = f(x)$ or $f(z) = f(y)$ whenever $z$ agrees with $x$ or $y$ respectively on $W$. We claim that $\Stab(\pi) \subseteq W^{-1}W$, and hence it is finite. Indeed, take any group element $\gamma \not \in W^{-1}W$. Then $W \cap W\gamma = \0$, so there is $z \in k^\G$ with $\rest{z}{W} = \rest{x}{W}$ and $\rest{(\gamma \cdot z)}{W} = \rest{y}{W}$. Since the set of all such points $z$ is open and $\Free(k^\G)$ is dense in $k^\G$, we may assume that $z \in \Free(k^\G)$. Then $f(z) = f(x) \neq f(y) = f(\gamma \cdot z)$, so $\gamma \not \in \Stab(\pi(z))$.
	\end{scproof}
	
	\begin{propcopy}{prop:Stab_conv}
	    For every finite normal subgroup $H \lhd \G$ and every $k \geq 2$, there is a continuous $\G$-equivariant map $\pi \colon \Free(k^\G) \to \mathcal{X}_H$ such that $\Stab(\pi) = H$.
	\end{propcopy}
	\begin{scproof}
	    Recall that $\mathcal{X}_H \subseteq 2^\G$ is the set of all functions $x \colon \G \to 2$ that are constant on the cosets of $H$. For $x \in \Free(k^\G)$, set $f(x)\defeq 0$ if $x(h) = 0$ for some $h \in H$ and $f(x) \defeq 1$ otherwise. Let $\pi \defeq \pi_f$ be the coding map for $f$. Then $\pi(x) \in \mathcal{X}_H$ and $H \leq \Stab(\pi(x))$ for all $x \in \Free(k^\G)$ by construction. Now take any group element $\gamma \not \in H$. There is a point $x \in k^\G$ with $x(h) = 0$ and $x(h\gamma) = 1$ for all $h \in H$. Moreover, since the set of all such points $x$ is open and $\Free(k^\G)$ is dense in $k^\G$, we may assume that $x \in \Free(k^\G)$. Then $f(x) = 0$ and $f(\gamma \cdot x) = 1$, which implies that $\gamma \not \in \Stab(\pi(x))$. %Fix any group element $\sigma \not \in \G$. Since $H$ is a subgroup of $\G$, we can consider the induced shift action $H \acts \Free(k^\G)$. Since $H$ is finite, the quotient space $X \defeq \Free(k^\G)/H$ is Polish. Let $\phi \colon \Free(k^\G) \to X$ denote the quotient map. The quotient group $\G/H$ acts on $X$ in the obvious way, and we let $G \defeq G(X, \set{H\sigma})$ be the Schreier graph of this action $\G/H \acts X$. Note that the maximum degree of $G$ is at most $2$. By Lemma~\ref{lemma:coloring} \ep{applied with $J_0 = J_1 = J_2 = \0$}, there is a partition $X = I_0 \sqcup I_1 \sqcup I_2$ into clopen sets that are independent in $G$.
	\end{scproof}
	
	\subsection{Maps to $\Col(F, \ell)$ with trivial $\Stab$}\label{subsec:col}
	
	In the introduction, we proved Corollary~\ref{corl:coloring} assuming that for all $\ell \geq |F|+1$, there is a continuous $\G$-equivariant map $\pi \colon \Free(3^\G) \to \Col(F,\ell)$ such that $\Stab(\pi) = \set{\mathbf{1}}$. Here we construct such $\pi$.
	
	\begin{lemma}
        Let $F \subset \G$ be a nonempty finite symmetric set with $\mathbf{1} \not \in F$. If $\ell \geq |F| + 1$, then there exists a continuous $\G$-invariant map $\pi \colon \Free(3^\G) \to \Col(F,\ell)$ such that $\Stab(\pi) = \set{\mathbf{1}}$.
	\end{lemma}
	\begin{scproof}
	    Let $G \defeq G(\Free(3^\G), F)$ be the Schreier graph of $\Free(3^\G)$. The maximum degree of $G$ is $|F| < \ell$. For $i \in \set{0,1}$, define
	    \[
	        J_i \,\defeq\, \set{x \in \Free(3^\G) \,:\, x(\mathbf{1}) = i \text{ and } x(\sigma) = 2 \text{ for all } \sigma \in F}.
	    \]
	    %$i < 3$, let $U_i \defeq \set{x \in \Free(3^\G) \,:\, x(\mathbf{1}) = i}$, and define
	    %\[
	    %    J_0 \,\defeq\, U_0 \cap (F \cdot U_2) \quad \text{and} \quad J_1 \,\defeq\, U_1 \cap (F \cdot U_2).
	    %\]
	    Note that the sets $J_0$, $J_1$ are clopen and independent in $G$. Let $g \colon J_0 \sqcup J_1 \to \ell$ be the map that sends every point in $J_0$ to $0$ and every point in $J_1$ to $1$. By Lemma~\ref{lemma:coloring}, $g$ can be extended to a continuous proper $\ell$-coloring $f \colon \Free(3^\G) \to \ell$ of $G$. Let $\pi \defeq \pi_f$ be the coding map for $f$. Then $\pi$ is a continuous $\G$-equivariant map from $\Free(3^\G)$ to $\Col(F, \ell)$, so it remains
	    %We need
	    to check that $\Stab(\pi) = \set{\mathbf{1}}$. Since $f(x) \neq f(\sigma \cdot x)$ for all $x \in \Free(3^\G)$ and $\sigma \in F$, we have $\Stab(\pi) \cap F = \0$. Now take any non-identity group element $\gamma \not \in F$. There is a point $x \in 3^\G$ with $x(\mathbf{1}) = 0$, $x(\gamma) = 1$, and $x(\sigma) = x(\sigma\gamma) = 2$ for all $\sigma \in F$. Since %the set of all such points is open and
	    $\Free(3^\G)$ is dense in $3^\G$, we may assume that $x \in \Free(3^\G)$. Then $x \in J_0$ and $\gamma \cdot x \in J_1$, so $f(x) = 0$ and $f(\gamma \cdot x) = 1$, which implies that $\gamma \not \in \Stab(\pi(x))$.
	\end{scproof}
	
	\section{A probabilistic tool: continuous Lov\'asz Local Lemma}\label{sec:LLL}
	
	The \emph{Lov\'asz Local Lemma}, or the \emph{LLL} for short, is a powerful tool in probabilistic combinatorics introduced by Erd\H{o}s and Lov\'asz \cite{EL}. Roughly speaking, the LLL asserts the existence of an object satisfying a set of constraints provided that \begin{enumerate*}[label=\ep{\emph{\alph*}}]
	    \item each constraint is unlikely to be violated by a \emph{random} object, and
	    \item the dependencies between the constraints are ``limited.''
	\end{enumerate*}
	The reader is referred to the books \cite{AS} by Alon and Spencer and \cite{MR} by Molloy and Reed for background on the LLL and a sample of its applications. 
	
	A recent trend involves establishing analogs of the LLL that are, in various senses, ``constructive'' \ep{as opposed to pure existence results}. This line of research was spurred by the breakthrough work of Moser and Tardos \cite{MT}, who developed an \emph{algorithmic} version of the LLL. A few other notable examples are the \emph{computable} LLL of Rumyantsev and Shen \cite{RSh}, the \emph{distributed} LLL of Fischer and Ghaffari \cite{FG}, the \emph{Borel} LLL of Cs\'{o}ka, Grabowski, M\'{a}th\'{e}, Pikhurko, and Tyros \cite{CGMPT}, and the \emph{measurable} LLL developed by the author \cite{Ber_LLL_shift, Ber_dist}. Here we will use the \emph{continuous} analog of the LLL introduced in \cite{Ber_cont}.
	
	Let $X$ be a set and let $C$ be a nonempty finite set. We refer to the elements of $C$ as \emphd{colors}. By a \emphd{$C$-coloring} of $X$ we mean a function $f \colon X \to C$. For a finite set $D \subseteq X$, an \emphd{$(X,C)$-constraint} \ep{or simply a \emphd{constraint} if $X$ and $C$ are clear} with \emphd{domain} $\dom(B) = D$ is a subset $B \subseteq C^D$. A $C$-coloring $f \colon X \to C$ \emphd{violates} a constraint $B$ with domain $D$ if $\rest{f}{D} \in B$, and \emphd{satisfies} $B$ otherwise. A \emphd{constraint satisfaction problem} \ep{a \emphd{CSP} for short} $\B$ on $X$ with range $C$, in symbols $\B \colon X \to^? C$, is an arbitrary set of $(X,C)$-constraints. A \emphd{solution} to a CSP $\B \colon X \to^? C$ is a $C$-coloring $f \colon X \to C$ that satisfies every constraint $B \in \B$. In other words, each constraint $B \in \B$ in a CSP $\B \colon X \to^? C$ is interpreted as a set of finite ``forbidden patterns'' that have to be avoided in a solution $f$.
	
	To each CSP $\B \colon X \to^? C$, we associate four numerical parameters:
	\[
	    \pr(\B), \qquad \de(\B), \qquad \vdeg(\B), \qquad \text{and} \qquad \ord(\B),
	\]
	defined as follows. For each $B \in \B$, the \emphd{probability} of $B$ is the quantity $\P[B] \defeq |B|/|C|^{|\dom(B)|}$. In other words, $\P[B]$ is the probability that $B$ is violated by a uniformly random $C$-coloring. Define
	\[
	    \pr(\B) \,\defeq\, \sup_{B \in \B} \P[B].
	\]%$\P[B]$ given by
	%\[
	%	\P[B] \,\defeq\, \frac{|B|}{k^{|\dom(B)|}} \,=\, \text{the probability that $B$ is violated by uniformly random $f \colon X \to k$}.
	%\]
	%Define $\pr(\B) \defeq \sup_{B \in \B} \P[B]$.
	The \emphd{maximum dependency degree} $\de(\B)$ of $\B$ is
	\[
	    \de(\B) \,\defeq\, \sup_{B \in \B} |\set{B' \in \B \,:\, B \neq B' \text{ and } \dom(B') \cap \dom(B) \neq \0}|.
	\]
	The \emphd{maximum vertex-degree} $\vdeg(\B)$ of $\B$ is
	\[
		\vdeg(\B) \,\defeq\, \sup_{x \in X} |\set{B\in \B \,:\, x \in \dom(B)}|,
	\]
	Finally, the \emphd{order} $\ord(\B)$ of $\B$ is \[\ord(\B) \,\defeq\, \sup_{B \in \B} |\dom(B)|.\]
	The LLL gives a condition, in terms of $\pr(\B)$ and $\de(\B)$, that guarantees that $\B$ has a solution:
	
	\begin{theo}[{{\textls{Lov\'asz Local Lemma}} \cite[Corollary 5.1.2]{AS}}]\label{theo:LLL}
		If $\B$ is a CSP such that
		\begin{equation}\label{eq:LLL}
			e \cdot \pr(\B) \cdot (\de(\B) + 1) \,\leq\, 1,
		\end{equation}
		where $e = 2.71\ldots$ is the base of the natural logarithm, then $\B$ has a solution.
	\end{theo}
	
	%Due to its origin in combinatorics, the LLL is often stated in the case when $\B$ is finite. However, a straightforward compactness argument shows that Theorem~\ref{theo:LLL} holds for infinite $\B$ as well \ep{see, e.g., \cite[proof of Theorem 5.2.2]{AS}}.
	
	The continuous LLL provides a condition similar to \eqref{eq:LLL} that guarantees the existence of a \emph{continuous} solution. To state it, we first need to define \emph{continuous CSPs}. Let $X$ be a zero-dimensional Polish space. For a set $B$ of functions $n \to C$ and a sequence of distinct points $x_0$, \ldots, $x_{n-1} \in X$, let $B(x_0, \ldots, x_{n-1})$ be the $(X,C)$-constraint given by
	$
	    B(x_0, \ldots, x_{n-1}) \,\defeq\, \set{\phi \circ \iota \,:\, \phi \in B},
	$
	where $\iota \colon \set{x_0, \ldots, x_{n-1}} \to n$ is the mapping $x_i \mapsto i$. A CSP $\B \colon X \to^? C$ is \emphd{continuous} if for all $n \in \N$, for every set $B$ of functions $n \to C$, and for all clopen subsets $U_1$, \ldots, $U_{n-1} \subseteq X$, the set %following set % is clopen:
	\[
		\left\{x_0 \in X \,:\, \begin{array}{c}\text{there are $x_1 \in U_1$, \ldots, $x_{n-1} \in U_{n-1}$ such that} \\ \text{$x_0$, \ldots, $x_{n-1}$ are distinct and $B(x_0, \ldots, x_{n-1}) \in \B$}\end{array}\right\}
	\]
	is also clopen. 
	\ep{This is analogous to the definition of a continuous graph from \S\ref{subsec:cont_graph}.}
	
	\begin{theo}[{\textls{Continuous LLL} \cite[Theorem 1.6]{Ber_cont}}]\label{theo:contLLL}
		Let $\B \colon X \to^? C$ be a continuous CSP on a zero-dimensional Polish space $X$. If
		\[%\begin{equation}\label{eq:bound}
			\pr(\B) \cdot \vdeg(\B)^{\ord(\B)} \,<\, 1,
		\]%\end{equation}
		then $\B$ has a continuous solution $f \colon X \to C$.
	\end{theo}
	
	Although the bound on $\pr(\B)$ required by Theorem~\ref{theo:contLLL} is much stronger than the one in the \hyperref[theo:LLL]{ordinary LLL}, it cannot be relaxed in general \ep{see the discussion after Theorem 1.6 in \cite{Ber_cont}}. Theorem~\ref{theo:contLLL} is particularly useful when $\ord(\B)$ is much smaller than $\vdeg(\B)$.
	
	\begin{remk}\label{remk:lists}
	In our proof of Theorem~\ref{theo:main}, we will need to use the \hyperref[theo:contLLL]{continuous LLL} in a somewhat more general setting than the one we have just described. Namely, suppose that in addition to a CSP $\B \colon X \to^? C$, we are given a mapping $x \mapsto C_x$ that assigns to each point $x \in X$ a nonempty subset $C_x \subseteq C$. Then we may seek a solution $f \colon X \to C$ to $\B$ with the extra property that $f(x) \in C_x$ for all $x \in X$. It then makes sense to consider colorings $f$ where each $f(x)$ is chosen uniformly at random from $C_x$, so the probability of a constraint $B \in \B$ should be defined by
	\[
	    \P[B] \,\defeq\, \frac{|B|}{\prod_{x \in \dom(B)} |C_x|}.
	\]
	Thankfully, there is a straightforward way to reduce this more general set-up to the one where all sets $C_x$ are equal to each other. Namely, let $n$ be a positive integer that is divisible by $1$, $2$, \ldots, $|C|$. For each nonempty subset $S \subseteq C$, fix an arbitrary function $h_S \colon n \to S$ such that the $h_S$-preimage of every element $s \in S$ has size exactly $n/|S|$. Then, instead of picking $f(x)$ uniformly at random from $C_x$, we can pick a uniformly random integer $0 \leq \phi(x) < n$ and set $f(x) \defeq h_{C_x}(\phi(x))$. For every constraint $B \in \B$, there is a corresponding $(X, n)$-constraint $B'$ given by
	\[
	    B' \,\defeq\, \set{\phi \colon \dom(B) \to n \,:\, \text{the map $x \mapsto h_{C_x}(\phi(x))$ is in $B$}}.
	\]
	In this way, we replace $\B$ by an ``equivalent'' CSP $\B' \defeq \set{B' \,:\, B \in \B} \colon X \to^? n$. Notice also that this construction preserves continuity, in the sense that if the CSP $\B$ and the function $x \mapsto C_x$ are continuous, then the CSP $\B'$ is also continuous, and, on the other hand, if $\phi \colon X \to n$ is a continuous solution to $\B'$, then $x \mapsto h_{C_x}(\phi(x))$ is a continuous solution to $\B$.
	\end{remk}
	
	\section{Proofs of Lemmas~\ref{lemma:cont+Stab} and \ref{lemma:col+Stab}}\label{sec:proof}
	
	\subsection{Inductive step lemma}\label{subsec:technical}
    
	In this subsection we prove a key technical result---Lemma~\ref{lemma:technical}---that will then be applied iteratively to obtain Lemmas~\ref{lemma:cont+Stab} and \ref{lemma:col+Stab}. It is a generalization of \cite[Lemma~4.1]{Ber_cont} and is proved using a similar strategy. However, its setting and proof are considerably more complicated than those of \cite[Lemma~4.1]{Ber_cont}, so the reader may find it instructive to look up the proof of \cite[Lemma~4.1]{Ber_cont} first (which, in contrast to the proof of Lemma~\ref{lemma:technical} in this paper, is less than two pages long). 
	
	\subsubsection{Informal overview}\label{subsubsec:informal}
	
	Since the statement of Lemma~\ref{lemma:technical} is technical and involves several changes of quantifiers, it would perhaps be helpful to first informally explain the intuition behind it.
	
	Let us start by describing the simplified setting of \cite[Lemma~4.1]{Ber_cont}. Suppose we are given a free continuous action $\G \acts X$ of $\G$ on a zero-dimensional Polish space $X$, and we want to construct a continuous mapping $f \colon X \to \ell$ for some $\ell \geq 2$ such that $\pi_f(X) \subseteq \Free(\ell^\G)$ \ep{where $\pi_{f}$ is the {coding map} for $f$, as defined in \S\ref{subsec:code}}.\footnote{Typically we would also want the closure of $\pi_f(X)$ to be contained in $\Free(\ell^\G)$, but for simplicity we will ignore this issue in the informal overview.} We may then proceed inductively, working with one non-identity group element at a time. At the stage when we consider an element $\gamma \in \G \setminus \set{\mathbf{1}}$, we will have already defined a continuous map $f_0 \colon C_0 \to \ell$ for some clopen subset $C_0 \subseteq X$. We say that the set $C_0$ has already been \emph{colored}. We will now extend $f_0$ to a continuous function $f \colon C_0 \sqcup C \to \ell$, where $C$ is some clopen subset of $X \setminus C_0$ (this is the set that we need to \emph{color} at this stage). The remaining \emph{uncolored} set $X \setminus (C_0 \sqcup C)$ is denoted by $U$. Our goal for this stage is to ensure that for every function $f' \colon X \to \ell$ extending $f$ and for all $x \in X$, we have $\gamma \not \in \Stab(\pi_{f'}(x))$.
	
	Since at this stage we only decide on the colors of the elements in $C$, we need the set $C$ to be ``large,'' while the remaining uncolored set $U$ should be ``small.'' The notion of ``largeness'' that we use here is \emph{syndeticity}: a set $A \subseteq X$ is called \emphd{$R$-syndetic} for a finite set $R \subset \G$ if $R^{-1} \cdot A = X$, i.e., if for each $x \in X$, there is some $\sigma \in R$ such that $\sigma \cdot x \in A$. We also say that a set $A \subseteq X$ is \emphd{$S$-separated} for a finite set $S \subset \G$ if $A$ is independent in the Schreier graph $G(X,S)$. Being $S$-separated for a suitable finite set $S \subset \G$ is our notion of ``smallness.'' With this set-up, \cite[Lemma~4.1]{Ber_cont} says that for every finite set $R \subset \G$, there is a \ep{much larger} finite set $S \subset \G$ such that a function $f$ with the desired properties exists provided that $C$ is $R$-syndetic and $U$ is $S$-separated. An important note to make here is that an $S$-separated set can still be $R'$-syndetic for some (even larger) finite set $R' \subset \G$; this is what allows the inductive construction to continue.
	
	Lemma~\ref{lemma:technical} in this paper generalizes the above result in a number of ways. Below we highlight the main differences between it and \cite[Lemma~4.1]{Ber_cont}.
	\begin{itemize}[wide]
	    \item The first difference is that we want the function $f \colon X \to \ell$ that we eventually construct to be a proper $\ell$-coloring of the Schreier graph $G(X,F)$ for some fixed finite set $F \subset \G$. Thankfully, Lemma~\ref{lemma:technical} assumes that $\ell$ (the number of colors) is strictly greater than the maximum degree of $G(X,F)$. By Lemma~\ref{lemma:coloring}, this implies that any partial proper $\ell$-coloring of $G(X,F)$ can be extended to a proper $\ell$-coloring of the entire graph, and so we can still use the inductive approach of building the coloring in stages.
	    
	    \item Second, our goal is no longer to make sure that $\pi_f(X) \subseteq \Free(\ell^\G)$. Instead, we are given a non-constant continuous function $\rho \colon \Col(F, \ell) \to m$. When $\Stab(\pi_\rho) = \set{\mathbf{1}}$, our goal becomes to find a continuous proper $\ell$-coloring $f \colon X \to \ell$ of $G(X,F)$ such that $(\pi_\rho \circ \pi_f)(X) \subseteq \Free(m^\G)$. This means that we now need to somehow control the function $\rho \circ \pi_f$ rather than $f$ itself. A key observation here is that, since $\rho$ is continuous and $\Col(F, \ell)$ is compact, the value $\rho \circ \pi_f$ only depends on the values of $f$ at finitely many points in the $\G$-orbit of $x$. %To this end, %notice that since $\rho$ is continuous and $\Col(F, \ell)$ is compact, there is some finite set $W \subset \G$ such that for every $x \in \Col(F, \ell)$, the value $\rho(x)$ is determined by the restriction of $x$ to $W$. This means that the value $(\rho \circ \pi_f)(x)$ for $x \in X$ is determined by the restriction of $f$ to the set $W \cdot x$. I
	    In view of this observation, it turns out to be useful to have a ``large'' set of points $x \in X$ with $D \cdot x \subseteq C$, where $D \subset \G$ is a certain finite set depending on $\rho$ (in particular, $D$ will remain fixed throughout the inductive construction). Indeed, if $D \cdot x \subseteq C$, then no point in $D \cdot x$ is colored yet, but all of them will become colored at the current stage of the construction. By choosing the coloring of $D \cdot x$ appropriately, we may then be able to control the value $(\rho \circ \pi_f)(x)$. %$W \cdot x \subseteq C$, since, by defining $f$ on $W \cdot x$ appropriately, we may then force $(\rho \cdot \pi_f)(x)$ to take a specific value. Actually, for various technical reasons, we need to have a ``large'' set of points $x \in X$ such that $D \cdot x \subseteq C$, where $D \subset \G$ is some other finite set that contains $W$ but is considerably larger. However, the set $D$ will still only depend on $\ell$, $F$, and $\rho$, and so it will remain unchanged throughout the inductive construction.
	    
	    Formally, we call a set $A \subseteq X$ is \emphd{$(R,D)$-syndetic} for finite $R$, $D \subset \G$ if the set $\set{x \in X \,:\, D \cdot x \subseteq A}$ is $R$-syndetic. Our ``largeness'' assumption now is that the set $C$ is $(R,D)$-syndetic for a certain set $D$ (depending on $\rho$) and some finite set $R$.
	    
	    \item Changing the ``largeness'' assumption on $C$ from `$R$-syndetic' to `$(R,D)$-syndetic' necessitates a change in the ``smallness'' assumption on $U$ as well. This is because, in order to be able to continue the inductive construction, we must ensure that the uncolored set $U$ is $(R',D)$-syndetic for some finite set $R' \subset \G$ (remember that the set $D$ will remain unchanged throughout the construction). Unfortunately, if $U$ is $S$-separated, it may be impossible to have even one element $x \in X$ such that $D \cdot x \subseteq U$, let alone an $R'$-syndetic set of such elements. To remedy this problem, we introduce a weaker notion of ``smallness'' that depends on the choice of the set $D$. Namely, we say that a set $A \subseteq X$ is \emphd{$(S,D)$-separated} for finite $S$, $D \subset \G$ if $A \subseteq D \cdot A'$ for some $S$-separated set $A' \subseteq X$. In the statement of Lemma~\ref{lemma:technical} we assume that $C$ is $(R, D)$-syndetic and $U$ is $(S,D)$-separated for some set $S$ that may depend on $R$. Crucially, being $(S,D)$-separated does not preclude being $(R',D)$-syndetic for some other finite set $R'$, which will help us continue the inductive process.
	    
	    \item The above discussion proceeded under the assumption that $\Stab(\pi_\rho) = \set{\mathbf{1}}$. In general, our goal is to find a continuous proper $\ell$-coloring $f \colon X \to \ell$ of $G(X,F)$ such that for all $x \in X$, we have $\Stab((\pi_\rho \circ \pi_f)(x)) = \Stab(\pi_\rho)$. Thankfully, since we are working with just one group element $\gamma \in \G \setminus \Stab(\pi_\rho)$ at a time, this results in only fairly minor technical changes to the proof.
	    
	    \item Finally, Lemma~\ref{lemma:technical} applies to actions $\G \acts X$ that are not necessarily free but only $S$-free for a large enough finite set $S \subset \G$ (where an action $\G \acts X$ is {$S$-free} if $S \cap \Stab(x) \subseteq \set{\mathbf{1}}$ for all $x \in X$). This feature of Lemma~\ref{lemma:technical} will be important in the proof of Lemma~\ref{lemma:cont+Stab}.
	\end{itemize}
	
	\subsubsection{The statement of the main technical lemma}
	
	Let us now proceed to the formal statement of Lemma~\ref{lemma:technical}. To begin with, we need a few definitions (some of which were already mentioned in \S\ref{subsubsec:informal}). Let $\G \acts X$ be a continuous action of $\G$ on a zero-dimensional Polish space. Given finite sets $R$, $S$, $D \subset \G$, we say that a subset $A \subseteq X$ is:
	\begin{itemize}
	    \item \emphd{$R$-syndetic} if $R^{-1} \cdot A = X$;
	    
	    \item \emphd{$(R,D)$-syndetic} if the set $\set{x \in X \,:\, D \cdot x \subseteq A}$ is $R$-syndetic;
	    
	    \item \emphd{$S$-separated} if it is independent in the Schreier graph $G(X,S)$; %for all distinct $x$, $y \in A$, $y \not \in S \cdot x$; % \ep{or, equivalently, if $A$ is independent in the Schreier graph $G(X,S)$};
	    
	    \item \emphd{$(S,D)$-separated} if $A \subseteq D \cdot A'$ for some $S$-separated set $A' \subseteq X$.
	\end{itemize}
	Let $F \subset \G$ be a finite symmetric subset with $\mathbf{1} \not \in F$ and let $\ell \geq |F| + 1$. Suppose that $f \colon X \pto \ell$ is a partial proper $\ell$-coloring of $G(X,F)$. The assumption $\ell \geq |F| + 1$ guarantees that $f$ can be extended to a proper $\ell$-coloring $f' \colon X \to \ell$ of the whole graph. \ep{Furthermore, by Lemma~\ref{lemma:coloring}, if $\dom(f)$ is clopen and $f$ is continuous, then $f'$ can be chosen to be continuous as well.} Given a function $\rho \colon \Col(F, \ell) \to m$, we say that two points $x$, $y \in X$ are \emphd{$(\rho, f)$-distinguished}, in symbols \[x \,\not \equiv_{\rho,f}\, y,\] if for every proper $\ell$-coloring $f' \colon X \to \ell$ that extends $f$, %we have
	\[
	    \rho(\pi_{f'}(x)) \,\neq\, \rho(\pi_{f'}(y)).
	\]
	Here $\pi_{f'}$ is the {coding map} for $f'$, as defined in \S\ref{subsec:code}. %Note that the assumption $\ell \geq |F| + 1$ guarantees that at least one such $f'$ exists.
	%point $x \in X$ is \emphd{$\rho$-decided by $f$} if for all proper $\ell$-colorings $f' \colon X \to \ell$ that extend $f$, the value $\rho(\pi_{f'}(x))$ is the same. \ep{Note that the assumption $\ell \geq |F| + 1$ guarantees that at least one such $f'$ exists.} The set of all points $x \in X$ that are $\rho$-decided by $f$ is denoted by $\dom(\rho, f)$. If $x \in \dom(\rho, f)$, then we let $\rho(x, f) \defeq \rho(\pi_{f'}(x))$ for any \ep{hence every} proper $\ell$-coloring $f' \colon X \to \ell$ extending $f$. Given a finite set $S \subset \G$, we say that two points $x$, $y \in X$ are \emphd{$(S,\rho,f)$-similar}, in symbols $x \equiv_{S, \rho, f} y$, if $\rho(\sigma \cdot x, f) = \rho(\sigma \cdot y, f)$ for all $\sigma \in S$ such that $\sigma \cdot x$, $\sigma \cdot y \in \dom(\rho, f)$. %Since the space $\Col(F, \ell)$ is compact, there exists a finite set $R \subset \G$ such that for every $x \in \Col(F, \ell)$, $\rho(x)$ is determined by the restriction of $x$ to $R$.
	Finally, recall that an action $\G \acts X$ is called \emphd{$S$-free} if $S \cap \Stab(x) \subseteq \set{\mathbf{1}}$ for all $x \in X$. We are now ready to state the lemma:
	
	\begin{lemma}\label{lemma:technical}
	    Fix the following data:
	    \begin{itemize}
	        \item a finite symmetric set $F \subset \G$ with $\mathbf{1} \not \in F$;
	        \item integers $\ell \geq |F| + 1$ and $m \geq 1$;
	        \item a non-constant continuous function $\rho \colon \Col(F, \ell) \to m$.
	    \end{itemize}
	    Then there exists a finite set $D \subset \G$ such that for every finite set $R \subset \G$ and every group element $\gamma \not \in \Stab(\pi_\rho)$, there is a finite set $S \subset \G$ with the following property.
	    
	    %\smallskip
	    
	    For every $S$-free continuous action $\G \acts X$ of $\G$ on a zero-dimensional Polish space $X$ and for every partition $X = C_0 \sqcup C \sqcup U$ of $X$ into clopen sets such that
	    \[
	        \text{$C$ is $(R,D)$-syndetic} \quad \text{and} \quad \text{$U$ is $(S,D)$-separated},
	    \]
	    setting $G \defeq G(X,F)$, we have that if $f_0 \colon C_0 \to \ell$ is a continuous proper $\ell$-coloring of $G[C_0]$, then $f_0$ can be extended to a continuous proper $\ell$-coloring $f \colon C_0 \sqcup C \to \ell$ of $G[C_0 \sqcup C]$ such that
	    \begin{equation}\label{eq:distinguish}
	        \forall x \in X\ \exists \sigma \in S\ : \quad \sigma \cdot x \not\equiv_{\rho, f} \sigma\gamma \cdot x.
	    \end{equation}
	    
	    \iffalse
	    Let $\G \acts X$ be an $S$-free continuous action of $\G$ on a zero-dimensional Polish space $X$. Set $G \defeq G(X, F)$. Let $X = C_0 \sqcup C \sqcup U$ be a partition of $X$ into clopen sets such that
	    \[
	        \text{$C$ is $(R,D)$-syndetic} \quad \text{and} \quad \text{$U$ is $(S,D)$-separated.}
	    \]
	    Suppose that $f_0 \colon C_0 \to \ell$ is a continuous proper $\ell$-coloring of $G[C_0]$. Then $f_0$ can be extended to a continuous proper $\ell$-coloring $f \colon C_0 \sqcup C \to \ell$ of $G[C_0 \sqcup C]$ such that
	    \begin{equation}\label{eq:distinguish}
	        \forall x \in X\ \exists \sigma \in S\ : \quad \sigma \cdot x \not\equiv_{\rho, f} \sigma\gamma \cdot x.
	    \end{equation}
	    \fi
	\end{lemma}
	
	%Recall that an action $\G \acts X$ is \emphd{$S$-free} if $S \cap \Stab(x) \subseteq \set{\mathbf{1}}$ for all $x \in X$.
	
	As discussed in \S\ref{subsubsec:informal}, in the notation of Lemma~\ref{lemma:technical}, the set $C_0$ is already \emph{colored}, the set $C$ is the one we need to \emph{color}, and the set $U$ is left \emph{uncolored}. 
	
	%In the remainder of this subsection, we prove Lemma~\ref{lemma:technical}.
	
	\subsubsection{Proof of Lemma~\ref{lemma:technical}}
	
	%\begin{scproof}
	    Let $F$, $\ell$, and $\rho$ be as in the statement of Lemma~\ref{lemma:technical}. For a finite set $W \subset \G$, let $\Col(F, \ell, W)$ denote the set of all proper $\ell$-colorings $W \to \ell$ of the \ep{finite} induced subgraph $G(\G, F)[W]$ of the Cayley graph $G(\G, F)$. Since the space $\Col(F, \ell)$ is compact and $\rho$ is continuous, there is a finite set $W \subset \G$ such that for every $x \in \Col(F, \ell)$, the value $\rho(x)$ is determined by the restriction of $x$ to $W$. In other words, there is a mapping $\tau \colon \Col(F, \ell, W) \to m$ such that
	    \[
	        \rho(x) \,=\, \tau(\rest{x}{W}) \quad \text{for all } x \in \Col(F, \ell).
	    \]
	    We shall assume, without loss of generality, that $W$ is symmetric and $\mathbf{1} \in W$.
	    
	    At this point, we can introduce a useful piece of notation. Let $\G \acts X$ be an action of $\G$ and let $f \colon X \pto \ell$ be  a partial proper $\ell$-coloring of $G(X,F)$. Suppose that $x \in X$ is a point such that $W \cdot x \subseteq \dom(f)$. Then for every proper $\ell$-coloring $f' \colon X \to \ell$ that extends $f$, the value $\rho(\pi_{f'}(x))$ is the same and determined by the restriction of $f$ to $W \cdot x$. We denote this value by $\rho_f(x)$. Given two points $x$, $y \in X$, we write $\rho_f(x) \neq^\ast \rho_f(y)$ if $W \cdot x$, $W \cdot y \subseteq \dom(f)$ and $\rho_f(x) \neq \rho_f(y)$. We also write $\rho_f(x) =^\ast \rho_f(y)$ if the statement $\rho_f(x) \neq^\ast \rho_f(y)$ is false. Clearly,
	    \[
	        \rho_f(x) \,\neq^\ast\, \rho_f(y) \quad \Longrightarrow \quad x \,\not\equiv_{\rho,f} \, y.
	    \]
	    The benefit of working with the relation $\rho_f(x) \neq^\ast \rho_f(y)$ instead of $x \not\equiv_{\rho,f} y$ is that it only depends on the values of $f$ on the finite set $W \cdot x \cup W \cdot y$.
	    
	    Next, we define a subgroup $H \leq \G$ as follows:
	    \[
	        H \,\defeq\, \set{h \in \G \,:\, \rho(x) = \rho(h \cdot x) \text{ for all } x \in \Col(F, \ell)}.
	    \]
	
	    \begin{claim}
	        We have $H \subseteq W^2 \cup WFW$. In particular, $H$ is finite.
	    \end{claim}
	    \begin{scproof}%\begin{stepproof}
	        The argument is similar to the proof of Proposition~\ref{prop:Stab}. Since $\rho$ is not constant, there are $x$, $y \in \Col(F, \ell)$ such that $\rho(x) \neq \rho(y)$. If $\gamma \not \in W^2 \cup WFW$, then the sets $W$ and  $W\gamma$ are disjoint and have no edges between them in the Cayley graph $G(\G,F)$, so we can define a proper partial coloring $W \sqcup W\gamma \to \ell$ that agrees with $x$ on $W$ and with $\gamma^{-1} \cdot y$ on $W\gamma$. Since $\ell \geq |F| + 1$, this partial coloring can be extended to a proper $\ell$-coloring $z \colon \G \to \ell$ of the Cayley graph $G(\G, F)$. Then $\rho(z) = \rho(x) \neq \rho(y) = \rho(\gamma \cdot z)$, and thus $\gamma \not \in H$.
	    \end{scproof}%\end{stepproof}
	
	%\begin{claim}\label{claim:normal}
	Note that %we have
	%\end{claim}
	%\begin{stepproof1}
	%    By definition,
	%\begin{align*}
	%        h \in \Stab(\pi_\rho) \,&\Longleftrightarrow \, \pi_\rho(x) = \pi_\rho(h \cdot x) \text{ for all } x \in \Col(F, \ell)\\
	%        &\Longleftrightarrow \, \rho(\gamma \cdot x) = \rho(\gamma h \cdot x) \text{ for all } x \in \Col(F, \ell) \text{ and } \gamma \in \G \\
	%       &\Longleftrightarrow \, \gamma h \gamma^{-1} \in H \text{ for all } \gamma \in \G,
	%\end{align*}
	%and hence
	$\Stab(\pi_\rho) = \set{h \,:\, \text{every conjugate of $h$ is in $H$}}$. Therefore, for each $h \in H \setminus \Stab(\pi_\rho)$, we can fix a group element $q(h)$ such that $q(h) h q(h)^{-1} \not \in H$. Let
	\[
	    Q \,\defeq\, \set{q(h) \,:\, h \in H \setminus \Stab(\pi_\rho)} \cup \set{\mathbf{1}}.
	\]
	We claim that Lemma~\ref{lemma:technical} holds with the following choice of $D$: %and let $D$ be the following finite symmetric subset of $\G$:
	\[
	    D \,\defeq\, (F \cup W \cup Q \cup Q^{-1})^{100}.
	\]
	%We will show that $D$ has the property described in Lemma~\ref{lemma:technical}.
	
	Let $R \subset \G$ be a finite set and let $\gamma \in \G \setminus \Stab(\pi_\rho)$. Upon replacing $R$ with a superset if necessary, %Without loss of generality,
	we may assume that $R$ is symmetric and $\mathbf{1} \in R$. Let $M$ be an arbitrary finite symmetric subset of $\G$ with $\mathbf{1} \in M$ whose size $|M|$ is sufficiently large as a function of $\ell$, $|D|$, and $|R|$ \ep{it will become clear from the rest of the proof how large $|M|$ needs to be}. %,  where $c > 0$ is so large that
	%\begin{equation}\label{eq:c}
	%	c > 2|D||Q||R|\text{\todo{???}}
	%\end{equation}
	Let \[N \,\defeq\, WDRM \cup MRDW,\] so $N$ is symmetric and contains $\mathbf{1}$. We shall prove Lemma~\ref{lemma:technical} for
	\[
	    S \,\defeq\, (N \cup \set{\gamma, \gamma^{-1}})^{1000}.
	\]
	
	Fix $X$, $G$, $C_0$, $C$, $U$, and $f_0$ as in Lemma~\ref{lemma:technical}. Since $C$ is $(R,D)$-syndetic, the set
	\[
	    C' \,\defeq\, \set{x \in X \,:\, D \cdot x \subseteq C}
	\]
	is $R$-syndetic. By Lemma~\ref{lemma:max}, we can find a clopen maximal $N^4$-separated subset $Z$ of $C'$. %We call the points $c \in C$ \emphd{centers}.
	The maximality of $Z$ means that $C' \subseteq N^4 \cdot Z$. Since $C'$ is $R$-syndetic, this implies that $Z$ is $N^4R$-syndetic. %For $f \colon C_0 \sqcup C \to 2$ and $x$, $y \in X$, we say that $x$ and $y$ are \emphd{$N$-similar} in $f$, in symbols $x \equiv^N_f y$, if
	%\[
	%\forall \eta \in N, \quad \set{\eta \cdot x,\, \eta \cdot y} \, \subseteq \, C_0 \sqcup C \quad \Longrightarrow \quad f(\eta \cdot x) \,=\, f(\eta \cdot y).
	%\]
	%For every $z \in Z$ and $\beta \in \G$, let
	%\[
	%    A(z, \beta) \,\defeq\, \set{\nu \in RM \,:\, (W \nu \cup W\nu \beta \gamma \beta^{-1}) \cdot z \subseteq C_0 \sqcup C}.
	%\]
	By Lemma~\ref{lemma:coloring}, we can extend $f_0$ to a continuous proper $\ell$-coloring \[
	    g \colon C_0 \sqcup (C \setminus (N \cdot Z)) \to \ell.
	\]
	Our goal now is to extend $g$ to a continuous proper $\ell$-coloring $f \colon C_0 \sqcup C \to \ell$ such that
	\begin{equation}\label{eq:Z}
	    \forall z \in Z \ \forall \beta \in N^5 \  \exists \nu \in DRM \,: \quad \rho_f(\nu \cdot z) \neq^\ast \rho_f(\nu\beta \gamma \beta^{-1} \cdot z).
	\end{equation}
	%where $\Delta \defeq \set{\beta \gamma \beta^{-1}\,:\, \beta \in N'R}$.
	
	\begin{claim}
	    If $f$ satisfies \eqref{eq:Z}, then it also satisfies \eqref{eq:distinguish}.
	\end{claim}
	\begin{scproof}%\begin{stepproof}
	    Suppose $f$ satisfies \eqref{eq:Z}. Take any $x \in X$. Since $Z$ is $N^4R$-syndetic, there is $\beta \in N^4R$ such that $\beta \cdot x \in Z$. Applying \eqref{eq:Z} with $z = \beta \cdot x$, we obtain $\nu \in DRM$ such that $\rho_f(\nu\beta \cdot x) \neq^\ast \rho_f(\nu \beta \gamma \cdot x)$. Since $\nu\beta\in S$, this yields %, i.e., there is $\eta \in N$ such that $\set{\eta \beta \cdot x, \eta \beta \gamma \cdot x} \subseteq \dom(f)$ and $f(\eta \beta \cdot x) \neq f(\eta \beta \gamma \cdot x)$. Since $\eta \beta \in N^5F = S$,
		the desired conclusion. %as $\nu\beta \in S$.
	\end{scproof}%\end{stepproof}
	
	What are the possible extensions of $g$ to $C_0 \sqcup C$? \ep{Here and in what follows, by an ``extension'' we mean an extension that is also a proper $\ell$-coloring.} Note that since $Z$ is $N^4$-separated, the sets $N \cdot z$ and $N \cdot z'$ for distinct $z$, $z' \in Z$ are disjoint and, moreover, there are no edges between them in $G$. This means that we may extend $g$ to each set $C \cap (N \cdot z)$ independently without creating any conflicts. Next, we consider the ways to extend $g$ to $C \cap (N \cdot z)$ for fixed $z \in Z$. To this end, let
	\[
	    N_z \,\defeq\, \set{\nu \in N \,:\, \nu \cdot z \in C}.
	\]
	Since the action $\G \acts X$ is $S$-free, there is a bijection $N_z \to C \cap (N \cdot z) \colon \nu \mapsto \nu \cdot z$. %Recall that $\Col(F, \ell, N_z)$ is the set of all proper $\ell$-colorings $N_z \to \ell$ of the \ep{finite} induced subgraph $G(\G, F)[N_z]$ of the Cayley graph $G(\G, F)$.
	Let $\pfun{N}{\ell}$ be the set of all partial functions $N \pto \ell$ and let $\Col(z) \subseteq \pfun{N}{\ell}$ be the subset comprising all functions $\phi \colon N_z \to \ell$ with the following properties:
	\begin{itemize}
	    \item $\phi$ is a proper coloring, i.e., $\phi(\nu) \neq \phi(\sigma \nu)$ whenever $\nu \in N_z$ and $\sigma \in F$ satisfy $\sigma \nu \in N_z$;
	    \item for all $\nu \in N_z$ and $\sigma \in F$, if $\sigma \nu \cdot z \in \dom(g)$, then $g(\sigma \nu \cdot z) \neq \phi(\nu)$.
	\end{itemize}
	In other words, $\phi \in \Col(z)$ if and only if we can extend $g$ to $C \cap (N \cdot z)$ by sending $\nu \cdot z$ to $\phi(\nu)$ for every $\nu \in N_z$. Thus, we can identify extensions of $g$ to $C_0 \sqcup C$ with functions \[Z \to \pfun{N}{\ell} \colon z \mapsto \phi_z\] such that $\phi_z \in \Col(z)$ for every $z \in Z$. Explicitly, given such a function $z \mapsto \phi_z$, we define the corresponding coloring $f \colon C_0 \sqcup C \to \ell$ extending $g$ via
	\[
	    f(\nu \cdot z) \,\defeq\, \phi_z(\nu) \quad \text{for all } z \in Z \text{ and } \nu \in N_z.
	\]
	Note that the assignment $z \mapsto \Col(z)$ is continuous. Also, if a function $z \mapsto \phi_z$ is continuous, then the corresponding coloring $f$ is continuous as well.
	
	Now we need to find a continuous function $z \mapsto \phi_z \in \Col(z)$ such that the corresponding coloring $f$ satisfies~\eqref{eq:Z}. We wish to phrase this problem as a CSP with domain $Z$. Consider any $z \in Z$ and $\beta \in N^5$. The truth of the statement
	\begin{equation}\label{eq:constraint}
	    \exists \nu \in DRM \ : \quad \rho_f(\nu \cdot z) \neq^\ast \rho_f(\nu \beta \gamma \beta^{-1} \cdot z)
	\end{equation}
	only depends on the restriction of $f$ to the set
	\[
	    (N \cup N\beta\gamma\beta^{-1}) \cdot z.
	\]
	The restriction of $f$ to $(N \cup N\beta\gamma\beta^{-1}) \cdot z$ is, in turn, determined by the functions $\phi_{z'}$ for all $z'$ in
	\[
	    \Delta(z, \beta) \,\defeq\, \set{z' \in Z \,:\, (N \cdot z') \cap (N \cup N\beta\gamma\beta^{-1}) \cdot z \neq \0}.
	\]
	Clearly, $z \in \Delta(z,\beta)$. Since $Z$ is $N^4$-separated, there can be at most one other $z' \in \Delta(z,\beta)$:
	
	\begin{claim}\label{claim:Delta}
	    For all $z \in Z$ and $\beta \in N^5$, $|\Delta(z, \beta)| \leq 2$.
	\end{claim}
	\begin{scproof}
	    By definition, $z' \in \Delta(z, \beta)$ if and only if
	    \[
	        z' \,\in\, Z \cap \left((N^2 \cup N^2 \beta \gamma \beta^{-1}) \cdot z\right).
	    \]
	    It remains to notice that each set of the form $N^2 \cdot x$, $x \in X$ includes at most one point from $Z$.
	\end{scproof}
	
	From the above discussion, it follows that for each $z \in Z$ and $\beta \in N^5$, we can form a constraint $B(z, \beta)$ with domain $\dom(B(z, \beta)) = \Delta(z, \beta)$ such that a mapping $z' \mapsto \phi_{z'}$ satisfies $B(z, \beta)$ if and only if \eqref{eq:constraint} holds for the corresponding coloring $f$. We then define a CSP $\B \colon Z \to^? \pfun{N}{\ell}$ by
	\[
	    \B \,\defeq\, \set{B(z,\beta) \,:\, z \in Z, \ \beta \in N^5}.
	\]
	Notice that the structure of the constraints involving any $z \in Z$ %, it is enough to know what happens on some finite subset of the $\G$-orbit of $z$; for instance, it suffices to know %the
	is fully determined by what happens on some finite subset of the $\G$-orbit of $z$. For instance, it is determined by %the following data are sufficient: 
	%\begin{enumerate*}[label=\ep{\itshape\alph*}] \item
	the intersections of $C$ and $Z$ with $S \cdot z$ and %\item
	the restriction of $g$ to $S \cdot z$%\end{enumerate*}
	. This implies that %some finite subset of the $\G$-orbit of $z$,
	the CSP $\B$ is continuous.
	
	We shall now apply the \hyperref[theo:contLLL]{continuous LLL} to argue that $\B$ has a continuous solution. To this end, we need to estimate $\pr(\B)$, $\ord(\B)$, and $\vdeg(\B)$. The latter two parameters are easy to bound:
	
	\begin{claim}\label{claim:ord}
	    $\ord(\B) \leq 2$.
	\end{claim}
	\begin{scproof}
	    Follows from Claim~\ref{claim:Delta}.
	\end{scproof}
	
	\begin{claim}\label{claim:vdeg}
	    $\vdeg(\B) \leq a |M|^{7}$, where $a > 0$ is a fixed function of $|D|$ and $|R|$.
	\end{claim}
	\begin{scproof}
	    For given $z \in Z$, we need to bound the number of pairs $(z', \beta)$ such that $z \in \Delta(z',\beta)$. Since $\beta \in N^5$, there are at most $|N|^5$ choices for $\beta$. Once $\beta$ is fixed, $z \in \Delta(z',\beta)$ implies
	    \[
	        z' \,\in\, (N^2 \cup \beta \gamma^{-1} \beta^{-1} N^2) \cdot z,
	    \]
	    so the number of choices for $z'$ is at most $2|N|^2$. Thus, the number of such pairs $(z', \beta)$ is at most
	    \[
	        |N|^5 \cdot 2|N|^2 \,=\, 2|N|^7 \,\leq\, 256|W|^7|D|^7|R|^{7}|M|^7 \,\leq\, 256|D|^{14}|R|^{7}|M|^{7},
	    \]
	    where we are using that $N = WDRM \cup MRDW$ and $D \supseteq W$.
	\end{scproof}
	
	It remains to bound $\pr(\B)$. Since for each $z \in Z$, the function $\phi_z$ must belong to $\Col(z)$, it makes sense to use the generalized setting for the LLL described in Remark~\ref{remk:lists}. That is, for the purposes of computing $\pr(\B)$, we consider picking each $\phi_z$ uniformly at random from $\Col(z)$.
	
	\begin{claim}\label{claim:probab}
	    $\pr(\B) \leq b^{|M|}$, where $0 < b < 1$ is a fixed function of $\ell$, $|D|$, and $|R|$.
	\end{claim}
	\begin{scproof}\stepcounter{ForClaims} \renewcommand{\theForClaims}{\ref{claim:probab}}
	    Take $z \in Z$ and $\beta \in N^5$. Define \[V \,\defeq\, C \cap (N \cdot \Delta(z,\beta)).\] We need to bound $\P[B(z, \beta)]$, i.e., the probability that the constraint $B(z, \beta)$ is violated when for each $z' \in \Delta(z, \beta)$, we pick a coloring $\phi_{z'} \in \Col(z')$ uniformly at random. If $\Delta(z,\beta) = \set{z}$, this means picking a uniformly random extension of $g$ to $N_z \cdot z = V$. If, on the other hand, $\Delta(z,\beta) = \set{z,z'}$ for some $z' \neq z$, then we uniformly at random extend $g$ to $N_z \cdot z$ and, independently, to $N_{z'} \cdot z'$. Since the sets $N_z \cdot z$ and $N_{z'} \cdot z'$ are disjoint and have no edges between them, this is the same as picking a uniformly random extension of $g$ to $V$. To summarize, regardless of whether $|\Delta(z,\beta)|$ is $1$ or $2$, our problem can be formulated as follows: Pick a uniformly random extension $\phi$ of $g$ to the set $V$. What is the probability that
	    \[
	        \forall \nu \in DRM \ : \quad \rho_\phi(\nu \cdot z) =^\ast \rho_\phi(\nu\beta\gamma\beta^{-1} \cdot z)?
	    \]
	    To bound this probability, we will focus on a carefully chosen subset of $DRM$.
	    
	    \begin{subclaim}\label{subclaim:E}
	        There exists a set $E \subseteq DRM$ with the following properties.
	        \begin{enumerate}[label=\ep{\itshape\alph*}]
	            \item\label{item:good} For all $\nu \in E$, $\nu\beta\gamma\beta^{-1}\nu^{-1} \not \in H$.
	            \item\label{item:uncolored} For all $\nu \in E$, $(F \cup W)^{10}\nu \cdot z \subseteq C$.
	            \item\label{item:colored} For all $\nu \in E$, $W\nu\beta\gamma\beta^{-1} \cdot z \subseteq C_0 \sqcup C$.
	            \item\label{item:disjoint} The sets $D\nu \cup W\nu\beta\gamma\beta^{-1}$ for $\nu \in E$ are pairwise disjoint.
	            \item\label{item:size} We have $|E| \geq  c |M|$, where $c>0$ is a fixed function of $|D|$ and $|R|$.
	        \end{enumerate}
	    \end{subclaim}
	    \begin{claimproof}[Proof of Subclaim~\ref{subclaim:E}]
	        Define the following sets:
	        \begin{align*}
	            E_1 \,&\defeq\, \set{\nu \in RM \,:\, \nu \cdot z \in C'};\\
	            E_2 \,&\defeq\, \set{\nu \in QE_1 \,:\, \nu\beta\gamma\beta^{-1}\nu^{-1} \not \in H};\\
	            E_3 \,&\defeq\, \set{\nu \in E_2 \,:\, W\nu\beta\gamma\beta^{-1} \cdot z \subseteq C_0 \sqcup C}.
	        \end{align*}
	        Let $E \subseteq E_3$ be a largest subset of $E_3$ such that the sets $D\nu \cup W\nu\beta\gamma\beta^{-1}$ for $\nu \in E$ are pairwise disjoint. We claim that $E$ is as desired. Conditions \ref{item:good}, \ref{item:colored}, and \ref{item:disjoint} hold by definition. To verify \ref{item:uncolored}, consider any $\nu \in E$. Since $\nu \in E_2$, we can write $\nu = q \nu'$ for some $q \in Q$ and $\nu' \in E_1$. By definition, $\nu' \cdot z \in C'$, which means that $D\nu' \cdot z \subseteq C$. Since $(F \cup W)^{10}Q \subseteq D$, this yields \ref{item:uncolored}.
	        
	        It remains to get a lower bound on $|E|$. Since the set $C'$ is $R$-syndetic, for each $\mu \in M$, we have $E_1 \cap R\mu \neq \0$. As every element of $E_1$ can belong to at most $|R|$ sets $R\mu$, we conclude that
	        \[
	            |E_1| \,\geq\, \frac{|M|}{|R|}.
	        \]
	        To obtain a lower bound on $|E_2|$, we observe that the set
	        \[
	            \Xi \,\defeq\, \set{\nu \in \G \,:\, \nu\beta\gamma\beta^{-1}\nu^{-1} \not \in H}
	        \]
	        is $Q$-syndetic. Indeed, take any $\nu \in \G$ and let $h \defeq \nu\beta\gamma\beta^{-1}\nu^{-1}$. Note that $h \not \in \Stab(\pi_\rho)$, because $h$ is a conjugate of $\gamma$ and $\gamma \not \in \Stab(\pi_\rho)$ by assumption. Hence, if $\nu$ itself is not in $\Xi$, then $h \in H \setminus \Stab(\pi_\rho)$, and hence we have an element $q(h) \in Q$ such that
	        \[
	            (q(h)\nu)\beta\gamma\beta^{-1}(q(h)\nu)^{-1}  \,=\, q(h) h q(h)^{-1} \,\not\in\, H,
	        \]
	        which means that $q(h)\nu \in \Xi$. Therefore, $E_2 \cap Q\nu \neq \0$ for all $\nu \in E_1$, so we can write
	        \[
	            |E_2| \,\geq\, \frac{|E_1|}{|Q|} \,\geq\, \frac{|E_1|}{|D|} \,\geq\, \frac{|M|}{|D||R|}. 
	        \]
	        To bound $|E_3|$, we use the fact that the set $U = X \setminus (C_0 \sqcup C)$ is $(S,D)$-separated. By definition, this means that there is an $S$-separated set $U'$ with $U \subseteq D \cdot U'$. The set $DWQRM\beta\gamma\beta^{-1} \cdot z$ can include at most one point of $U'$, so there are at most $|D||W| \leq |D|^2$ elements in $E_2 \setminus E_3$. Hence,
	        \[
	            |E_3| \,\geq\, |E_2| - |D|^2 \,\geq\, \frac{|M|}{|D||R|} - |D|^2 \,\geq\, \frac{|M|}{2|D||R|},
	        \]
	        where in the last inequality we use that $|M|$ is large enough in terms of $|D|$ and $|R|$. Finally, consider the graph with vertex set $E_3$ where two distinct elements $\nu$, $\nu' \in E_3$ are adjacent if and only if
	        \[
	            \big(D\nu \cup W\nu\beta\gamma\beta^{-1}\big) \, \cap \, \big(D\nu' \cup W\nu'\beta\gamma\beta^{-1}\big) \,\neq\, \0.
	        \]
	        The maximum degree of this graph is less than $4|D|^2$. Therefore, $|E|$---i.e., the largest size of an independent set in this graph---is at least
	        \[
	            \frac{|E_3|}{4|D|^2} \,\geq\, \frac{|M|}{8|D|^3|R|}. \qedhere
	        \]
	    \end{claimproof}
	    
	    Fix a set $E$ given by Subclaim~\ref{subclaim:E}. We will bound from above the probability that a uniformly random extension $\phi$ of $g$ to $V$ satisfies
	    \begin{equation}\label{eq:bad}
	        \forall \nu \in E \ : \quad \rho_\phi(\nu \cdot z) =^\ast \rho_\phi(\nu\beta\gamma\beta^{-1} \cdot z).  
	    \end{equation}
	    By \ref{item:uncolored} and \ref{item:colored}, for each $\nu \in E$, the sets $W\nu \cdot z$ and $W \nu\beta\gamma\beta^{-1} \cdot z$ are contained in the domain of $\phi$, so the values $\rho_\phi(\nu \cdot z)$ and $\rho_\phi(\nu\beta\gamma\beta^{-1} \cdot z)$ are defined. Thus, condition \eqref{eq:bad} is equivalent to
	    \begin{equation}\label{eq:bad1}
	        \forall \nu \in E \ : \quad \rho_\phi(\nu \cdot z) = \rho_\phi(\nu\beta\gamma\beta^{-1} \cdot z).  
	    \end{equation}
	    Split the set $V$ into two subsets:
	    \[
	        V_0 \,\defeq\, V \setminus \big((F \cup W)^{10}E \cdot z\big) \quad \text{and} \quad V_1 \,\defeq\, (F \cup W)^{10}E \cdot z.
	    \]
	    \ep{Note that, by \ref{item:uncolored}, $V_1$ is indeed a subset of $V$.} Let $\Ext_0$ denote the set of all extensions of $g$ to $V_0$. For each $\psi \in \Ext_0$, let $\Ext(\psi)$ be the set of all extensions of $\psi$ to $V_1$, and let $\Ext'(\psi) \subseteq \Ext(\psi)$ be the set of all such extensions $\phi$ that satisfy \eqref{eq:bad1}. Then the probability we wish to bound is equal to
	    \[
	        p \,\defeq\, \frac{\sum_{\psi \in \Ext_0} |\Ext'(\psi)|}{\sum_{\psi \in \Ext_0} |\Ext(\psi)|}.
	    \]
	    Fix any $\psi \in \Ext_0$. For each $\nu \in E$, let
	    \[
	        V_\nu \,\defeq\, (F \cup W)^{10}\nu \cdot z,
	    \]
	    and let $\Ext(\psi, \nu)$ be the set of all extensions of $\psi$ to $V_\nu$. By \ref{item:disjoint}, the sets $V_\nu$ for different $\nu \in E$ are disjoint and have no edges between them, so an extension of $\psi$ to $V_1$ is obtained by putting together arbitrary extensions of $\psi$ to each set $V_\nu$. Therefore,
	    \[
	        |\Ext(\psi)| \,=\, \prod_{\nu \in E} |\Ext(\psi, \nu)|.
	    \]
	    It follows from \ref{item:colored} and \ref{item:disjoint} that for every $\nu \in E$, the set $W \nu\beta\gamma\beta^{-1} \cdot z$ is included in $C_0 \sqcup V_0 \sqcup V_\nu$. Therefore, if $\xi \in \Ext(\psi, \nu)$, then both $\rho_\xi(\nu \cdot z)$ and $\rho_\xi(\nu\beta\gamma\beta^{-1} \cdot z)$ are defined. This allows us to define $\Ext'(\psi, \nu)$ as the set of all $\xi \in \Ext(\psi, \nu)$ such that
	    \[
	        \rho_\xi(\nu \cdot z) = \rho_\xi(\nu\beta\gamma\beta^{-1} \cdot z).  
	    \]
	    Then we have
	    \[
	        |\Ext'(\psi)| \,=\, \prod_{\nu \in E} |\Ext'(\psi, \nu)|.
	    \]
	    %The final step of the argument is the following subclaim.
	    
	    \begin{subclaim}\label{subclaim:exists}
	        For every $\psi \in \Ext_0$ and $\nu \in E$, $\Ext'(\psi, \nu) \neq \Ext(\psi, \nu)$.
	    \end{subclaim}
	    \begin{claimproof}[Proof of Subclaim~\ref{subclaim:exists}]
	        We have to show that for each $\psi \in \Ext_0$ and $\nu \in E$, there is an extension $\xi$ of $\psi$ to $V_\nu$ such that $\rho_\xi(\nu \cdot z) \neq \rho_\xi(\nu\beta\gamma\beta^{-1} \cdot z)$. We consider two cases.
	        
	        \smallskip
	        
	        \noindent \textbf{Case 1:} \textsl{The sets $(W \cup FW)\nu$ and $W\nu\beta\gamma\beta^{-1}$ are disjoint}. In this case we first take an arbitrary extension $\eta$ of $\psi$ to $W\nu\beta\gamma\beta^{-1} \cdot z$ and let $i \defeq \rho_\eta(\nu\beta\gamma\beta^{-1} \cdot z)$. Since the action $\G \acts X$ is $S$-free, the sets $(W \cup FW)\nu \cdot z$ and $W\nu\beta\gamma\beta^{-1} \cdot z$ are also disjoint, so the set $(W \cup FW)\nu \cdot z$ is entirely uncolored in $\eta$. Therefore, we may place an arbitrary proper $\ell$-coloring on $W\nu \cdot z$. In particular, since $\rho$ is not constant, we may color $W\nu \cdot z$ so that in the resulting coloring $\xi$, the value $\rho_\xi(\nu \cdot z)$ is distinct from $i$, as desired.
	        
	        \smallskip
	        
	        \noindent \textbf{Case 2:} \textsl{The sets $(W \cup FW)\nu$ and $W\nu\beta\gamma\beta^{-1}$ are not disjoint}. This implies that
	        \[
	            (W \cup FW)\nu \beta\gamma\beta^{-1} \,\subseteq \, V_\nu,
	        \]
	        which means that we can place an arbitrary proper $\ell$-coloring on the set
	        \[
	            \left(W\nu \cup W\nu \beta\gamma\beta^{-1}\right) \cdot z \,=\, \left(W \cup Wh\right) \cdot \nu z,
	        \]
	        where $h \defeq \nu\beta\gamma\beta^{-1}\nu^{-1}$. Since, by \ref{item:good}, $h \not \in H$, there is a coloring $x \in \Col(F,\ell)$ such that $\rho(x) \neq \rho(h \cdot x)$. It remains to copy the restriction of $x$ to $W \cup Wh$ onto the set $\left(W \cup Wh\right) \cdot \nu z$.
	    \end{claimproof}
	    
	    Now we are ready to derive the desired bound on $p$. Note that for all $\psi \in \Ext_0$ and $\nu \in E$,
	    \[
	        |\Ext(\psi, \nu)| \,\leq\, \ell^{|V_\nu|} \,\leq\, \ell^{|D|}.
	    \]
	    Therefore, by Subclaim~\ref{subclaim:exists},
	    \[
	        |\Ext'(\psi, \nu)| \,\leq\, |\Ext(\psi, \nu)| - 1 \,\leq\, \left(1 - \ell^{-|D|}\right)|\Ext(\psi, \nu)|.
	    \]
	    Thus,
	    \[
	        |\Ext'(\psi)| \,=\, \prod_{\nu \in E} |\Ext'(\psi, \nu)| \,\leq\, \left(1 - \ell^{-|D|}\right)^{|E|} \prod_{\nu \in E} |\Ext(\psi, \nu)| \,=\, \left(1 - \ell^{-|D|}\right)^{|E|} |\Ext(\psi)|,
	    \]
	    and hence, using \ref{item:size}, we get
	    \[
	        p \,=\, \frac{\sum_{\psi \in \Ext_0} |\Ext'(\psi)|}{\sum_{\psi \in \Ext_0} |\Ext(\psi)|} \,\leq\, \left(1 - \ell^{-|D|}\right)^{|E|} \,\leq\, \left(1 - \ell^{-|D|}\right)^{c|M|}. \qedhere
	    \]
	    %where the last inequality follows by \ref{item:size}.
	\end{scproof}
	
	The stage is now set for an application of the \hyperref[theo:contLLL]{continuous LLL}. Recall that our goal is to find a continuous solution to the CSP $\B$. By Theorem~\ref{theo:contLLL}, $\B$ has a continuous solution provided that
	\[
	    \pr(\B) \cdot \vdeg(\B)^{\ord(\B)}  \,<\, 1.
	\]
	Using Claims \ref{claim:ord}, \ref{claim:vdeg}, and \ref{claim:probab}, we may write
	\[
	    \pr(\B) \cdot \vdeg(\B)^{\ord(\B)} \,\leq\, b^{|M|} (a|M|^7)^2 \,=\, a^2|M|^{14}b^{|M|},
	\]
	where $a$ and $b$ are functions of $\ell$, $|D|$, and $|R|$ given by Claims \ref{claim:vdeg} and \ref{claim:probab}. It remains to observe that since $0 < b < 1$, the desired inequality
	\[
	    a^2|M|^{14}b^{|M|} \,<\, 1
	\]
	holds as long as $|M|$ is large enough (as a function of $\ell$, $|D|$, and $|R|$).
	
	%\end{scproof}
	
	\subsection{Proof of Lemma~\ref{lemma:col+Stab}}\label{subsec:col+Stab}
	
	In the remainder of \S\ref{sec:proof}, we prove Lemmas~\ref{lemma:cont+Stab} and \ref{lemma:col+Stab}. We start with Lemma~\ref{lemma:col+Stab}, since its proof follows more straightforwardly from Lemma~\ref{lemma:technical}. For convenience, we restate Lemma~\ref{lemma:col+Stab} here:
	
	\begin{lemmacopy}{lemma:col+Stab}
	       Let $F \subset \G$ be a finite symmetric set with $\mathbf{1} \not \in F$. Fix integers $\ell \geq |F| + 1$, $m \geq 1$ and let $\pi \colon \Col(F,\ell) \to m^\G$ be a continuous $\G$-equivariant map. Then there exists a subshift $\mathcal{Z} \subseteq  \Col(F,\ell)$ with the following properties:
	    \begin{itemize}
	        \item for all $z \in \mathcal{Z}$, $\Stab(\pi(z)) = \Stab(\pi)$; and
	        \item every free Borel action $\G \acts X$ on a Polish space admits a Borel $\G$-equivariant map $X \to \mathcal{Z}$.
	    \end{itemize}
	\end{lemmacopy}
	\begin{scproof}%\stepcounter{ForClaims} \renewcommand{\theForClaims}{\ref{lemma:col+Stab}}
	    %This is a modification of the proof of Theorem~\ref{theo:STD} given in \cite[\S4.B]{Ber_cont}, with Lemma~\ref{lemma:technical} replacing \cite[Lemma 4.1]{Ber_cont}.
	    \stepcounter{ForClaims} \renewcommand{\theForClaims}{\ref{lemma:col+Stab}}
	    We may assume that $\pi$ is not constant, since otherwise we can just take $\mathcal{Z} = \Col(F, \ell)$. Let $\rho \colon \Col(F, \ell) \to m$ be the \ep{continuous} function such that $\pi = \pi_\rho$. Explicitly, %we have
	    \[
	        \rho(x) \,\defeq\, \pi(x)(\mathbf{1}) \quad \text{for all } x \in \Col(F, \ell).
	    \]
	    Let $D \subset \G$ be the finite subset given by Lemma~\ref{lemma:technical} applied to $F$, $\ell$, $m$, and $\rho$. Without loss of generality, we may assume that $D$ is symmetric and contains $\mathbf{1}$.
	    
	    Fix an arbitrary enumeration $\gamma_0$, $\gamma_1$, \ldots{} of the elements of $\G \setminus \Stab(\pi)$. We recursively define a sequence of finite sets $T_0$, $R_0$, $S_0$, $T_1$, $R_1$, $S_1$, \ldots{} $\subset \G$ as follows. Set $T_0 \defeq \set{\mathbf{1}}$. Once $T_n$ is defined, let $Q_n \subset \G$ be any finite subset of $\G$ with $|Q_n| > |T_n||D|^2$ and set $R_n \defeq T_nQ_n$. Let $S_n$ be the set $S$ produced by Lemma~\ref{lemma:technical} applied with $R = R_n$ and $\gamma = \gamma_n$. After replacing $S_n$ with a superset if necessary, we arrange so that $S_n$ is symmetric, contains $\mathbf{1}$, and satisfies
	    \[
	        S_n \,\supseteq\, D^2R_nR_n^{-1}D^2. 
	    \]
	    Finally, we let $T_{n+1} \defeq S_nT_n$. This construction is done so that the following claim holds:
	
	\begin{smallclaim}\label{claim:splitD}
		Let $\G \acts X$ be a free continuous action of $\G$ on a zero-dimensional Polish space and let $V \subseteq X$ be a $(T_n,D)$-syndetic clopen set. Then there is a partition $V = C \sqcup U$ into two clopen subsets, where $C$ is $(R_n,D)$-syndetic, while $U$ is both $(S_n,D)$-separated and $(T_{n+1},D)$-syndetic.
	\end{smallclaim}
	\begin{claimproof}[Proof of Claim~\ref{claim:splitD}]
		Since $V$ is $(T_n, D)$-syndetic, the set
		\[V' \,\defeq\, \set{x \in X \,:\, D \cdot x \subseteq V}\]
		is $T_n$-syndetic. Note that $V'$ is clopen, so we can apply Lemma~\ref{lemma:max} to obtain a clopen maximal $S_n$-separated subset $U' \subseteq V'$. The maximality of $U'$ means that $V' \subseteq S_n \cdot U'$, and, since $V'$ is $T_n$-syndetic and $T_{n+1} = S_n T_n$, this implies that $U'$ is $T_{n+1}$-syndetic. Define
		\[
		    U \,\defeq\, D \cdot U' \quad \text{and} \quad C \,\defeq\, V \setminus U.
		\]
		Since $U'$ is $S_n$-separated and $T_{n+1}$-syndetic, $U$ is $(S_n, D)$-separated and $(T_{n+1},D)$-syndetic. %, as desired.
		It remains to verify that $C$ is $(R_n, D)$-syndetic. To this end, let
		\[
		    C' \,\defeq\, V' \setminus (D^2 \cdot U').
	    \]
	    Then the set $D \cdot C'$ is disjoint from $U$, and hence it is a subset of $C$. We claim that $C'$ is $R_n$-syndetic, which implies that $C$ is $(R_n, D)$-syndetic, as desired. Take any $x \in X$. We need to argue that $R_n \cdot x$ contains a point in $C'$. Recall that $R_n = T_n Q_n$. Since $V'$ is $T_n$-syndetic, $(T_n q \cdot x) \cap V' \neq \0$ for all $q \in Q_n$. As every element of $V'$ belongs to at most $|T_n|$ sets of the form $T_n q \cdot x$, we conclude that
		\begin{equation}\label{eq:lowerC'}
		    |(R_n \cdot x) \cap V'| \,\geq\, \frac{|Q_n|}{|T_n|} \,>\, |D|^2.
		\end{equation}
		On the other hand, since $U'$ is $S_n$-separated and $S_n \supseteq D^2R_nR_n^{-1}D^2$, there is at most one point $y \in U'$ with $(R_n \cdot x) \cap (D^2 \cdot y) \neq \0$. Therefore,
		\begin{equation}\label{eq:upperC'}
		    |(R_n \cdot x) \cap (D^2 \cdot U')| \,\leq\, |D|^2.
		\end{equation}
		From \eqref{eq:lowerC'} and \eqref{eq:upperC'}, it follows that $|(R_n \cdot x) \cap C'| > 0$, and we are done.
	\end{claimproof}
	
	For each $n \in \N$, let $\mathcal{Z}_n \subseteq \Col(F, \ell)$ be the set of all colorings $z \in \Col(F, \ell)$ such that
	\[
		\exists \sigma \in S_n \ : \quad \rho(\sigma \cdot z) \,\neq\, \rho(\sigma\gamma_n \cdot z).
	\]
	Since $\rho$ is continuous, $\mathcal{Z}_n$ is a relatively clopen set in $\Col(F, \ell)$. Also, if $z \in \mathcal{Z}_n$, then $\pi(z) \neq \pi(\gamma_n \cdot z)$, i.e., $\gamma_n \not \in \Stab(\pi(z))$. Therefore, the set
	\[
	    \mathcal{Z} \,\defeq\, \bigcap_{n = 0}^\infty \bigcap_{\delta \in \G} (\delta \cdot \mathcal{Z}_n)
	\]
	is a subshift contained in $\Col(F, \ell)$ with the property that $\Stab(\pi(z)) = \Stab(\pi)$ for every $z \in \mathcal{Z}$. %\ep{Notice that we have not yet shown that $\mathcal{Z}$ is nonempty.}
	To finish the proof of Lemma~\ref{lemma:col+Stab}, it remains to argue that every free Borel action of $\G$ on a Polish space admits a Borel $\G$-equivariant map to $\mathcal{Z}$.
	
	Let $\G \acts X$ be a free Borel action on a Polish space $X$. By replacing the topology on $X$ with a finer one if necessary, we may assume that $X$ is zero-dimensional and the action $\G \acts X$ is continuous \cite[\S13]{KechrisDST}. Set $G \defeq G(X,F)$. Iterative applications of Claim~\ref{claim:splitD} yield a sequence of clopen subsets $C_0$, $U_0$, $C_1$, $U_1$, $C_2$, $U_2$, \ldots{} of $X$ such that $C_0 = \0$, $U_0 = X$, and for all $n \in \N$,
	\begin{itemize}
		\item $U_n = C_{n+1} \sqcup U_{n+1}$, and
		\item the set $C_{n+1}$ is $(R_n, D)$-syndetic, while $U_{n+1}$ is $(S_n,D)$-separated and $(T_{n+1},D)$-syndetic.
	\end{itemize}
	%Let $C \defeq \bigsqcup_{n=0}^\infty C_n$. Note that the set $C$ is open.
	We then use Lemma~\ref{lemma:technical} repeatedly to obtain an increasing sequence $\0 = f_0 \subseteq f_1 \subseteq f_2 \subseteq \ldots$ such that for each $n \in \N$, $f_{n+1}$ is a continuous proper $\ell$-coloring of $G[C_0 \sqcup C_1 \sqcup \ldots \sqcup C_{n+1}]$ satisfying
	\begin{equation}\label{eq:not_sim_D}
	    \forall x \in X\ \exists \sigma \in S_n\ : \quad \sigma \cdot x \not\equiv_{\rho, f_{n+1}} \sigma\gamma_n \cdot x.
	\end{equation}
	Let $f \colon X \to \ell$ be any Borel proper $\ell$-coloring of $G$ extending $\acup_{n=0}^\infty f_n$. \ep{Such a coloring $f$ exists by the Borel version of Lemma~\ref{lemma:coloring}, the proof of which can be found, e.g., in \cite[Corollary 2.2]{BerConley}.} We claim that the map $\pi_f$ is as desired, i.e., $\pi_f(X) \subseteq \mathcal{Z}$. Indeed, since $f$ is a proper $\ell$-coloring of $G$, we have $\pi_f(X) \subseteq \Col(F, \ell)$. Since $f$ is an extension of $f_{n+1}$, \eqref{eq:not_sim_D} implies that for all $x \in X$,
	\[
		\exists \sigma \in S_n\ : \quad  \rho(\pi_f(\sigma \cdot x)) \,\neq\, \rho(\pi_f(\sigma\gamma_n \cdot x)).
	\]
	Therefore, $\pi_f(x) \in \mathcal{Z}_n$ for all $x \in X$ and $n \in \N$, i.e., $\pi_f(X) \subseteq \mathcal{Z}$, as desired.
	\end{scproof}

	\subsection{Proof of Lemma~\ref{lemma:cont+Stab}}\label{sec:cont+Stab}
	
	\subsubsection{Preparation}
	
	To prove Lemma~\ref{lemma:cont+Stab}, we will need the following consequence of Lemma~\ref{lemma:technical}:
	
	\begin{corl}\label{corl:less_technical}
	    Fix an integer $k \geq 2$, a finite set $R \subset \G$, and a group element $\gamma \neq \mathbf{1}$. Then there is a finite set $S \subset \G$ with the following property. Let $\G \acts X$ be an $S$-free continuous action of $\G$ on a zero-dimensional Polish space $X$ and let $X = C_0 \sqcup C \sqcup U$ be a partition of $X$ into clopen sets, where $C$ is $R$-syndetic and $U$ is $S$-separated. Then every continuous function $f_0 \colon C_0 \to k$ can be extended to a continuous function $f \colon C_0 \sqcup C \to k$ such that
	    \[%\begin{equation}\label{eq:distinguish}
	        \forall x \in X\ \exists \sigma \in S\ : \quad \set{\sigma \cdot x,\, \sigma \gamma \cdot x} \,\subseteq\, C_0 \sqcup C \text{ and } f(\sigma \cdot x) \,\neq\, f(\sigma\gamma \cdot x).
	    \]%\end{equation}
	\end{corl}
	\begin{scproof}
	    The Cayley graph $G(\G, \0)$ has no edges, so $\Col(\0, k) = k^\G$. Using this observation, we apply Lemma~\ref{lemma:technical} with $F = \0$, $\ell = m = k$, and the function $\rho \colon k^\G \to k$ given by \[\rho(x) \defeq x(\mathbf{1}) \quad \text{for all } x \in k^\G.\] In the notation of \S\ref{subsec:technical}, we can then take $W = \set{\mathbf{1}}$. It is also clear that $H = \set{\mathbf{1}}$, which implies that $Q = \set{\mathbf{1}}$, and hence $D = \set{\mathbf{1}}$ as well. The desired conclusion now follows by Lemma~\ref{lemma:technical}.
	\end{scproof}
	
	We remark that %using Lemma~\ref{lemma:technical} in this way to prove Lemma~\ref{lemma:less_technical} is, in some sense, an overkill. This is because
	Corollary~\ref{corl:less_technical} is a slight strengthening of \cite[Lemma 4.1]{Ber_cont} and can be established using essentially the same proof \ep{which is considerably less technical than the proof of Lemma~\ref{lemma:technical} given in \S\ref{subsec:technical}}. Specifically, \cite[Lemma~4.1]{Ber_cont} is the special case of Corollary~\ref{corl:less_technical} when $k = 2$ and the action $\G \acts X$ is assumed to be free rather than $S$-free. %The proof of \cite[Lemma 4.1]{Ber_cont} yields Lemma~\ref{lemma:less_technical} as well and is considerably simpler than the proof of Lemma~\ref{lemma:technical} given in \S\ref{subsec:technical}. However, since we have already established Lemma~\ref{lemma:technical}, we have no reason not to use it now.
	
	The following statement is proved by applying Corollary~\ref{corl:less_technical} iteratively:
	
	\begin{lemma}\label{lemma:extension}
	    Fix integers $k \geq 2$ and $m \geq 1$ and let $\rho \colon \Free(k^\G) \to m$ be a continuous function. Given finite subsets $D$, $T \subset \G$, there exist:
	    \begin{itemize}
	        \item a subshift $\mathcal{Z} \subseteq k^\G$,
	        \item a continuous map $\tilde{\rho} \colon \mathcal{Z} \to m$, and
	        \item a finite set $S \subset \G$
	    \end{itemize}
	    with the following properties.
	    
	    \begin{enumerate}[label=\ep{\normalfont\Roman*}]
	        \item\label{item:I} For each $z \in \mathcal{Z}$, there is $z^\ast \in \Free(k^\G)$ such that
	        \begin{enumerate}[label=\ep{\normalfont{}I\itshape\alph*}]
	            \item\label{item:Ia} for all $\delta \in D$, $z(\delta) = z^\ast(\delta)$; and
	            \item\label{item:Ib} for all $\delta \in D$, $\tilde{\rho}(\delta \cdot z) = \rho(\delta \cdot z^\ast)$.
	        \end{enumerate}
	        \item\label{item:II} Let $\G \acts X$ be an $S$-free continuous action of $\G$ on a zero-dimensional Polish space $X$ and let $C_0 \subseteq X$ be a clopen set such that its complement $X \setminus C_0$ is $T$-syndetic. Then every continuous map $f_0 \colon C_0 \to k$ can be extended to a continuous map $f \colon X \to k$ with $\pi_f(X) \subseteq \mathcal{Z}$.
	    \end{enumerate}
	\end{lemma}
	\begin{scproof}%[ of Claim~\ref{claim:extension}]
	\stepcounter{ForClaims} \renewcommand{\theForClaims}{\ref{lemma:extension}}
        This argument is a modification of the proof of Theorem~\ref{theo:cont} given in \cite[\S4.C]{Ber_cont} and is similar to the proof of Lemma~\ref{lemma:col+Stab} given in \S\ref{subsec:col+Stab}. Fix an enumeration $\gamma_0$, $\gamma_1$, \ldots{} of the non-identity elements of $\G$. We recursively define a sequence of finite sets $T_0$, $R_0$, $S_0$, $T_1$, $R_1$, $S_1$, \ldots{} $\subset \G$ as follows. Let $T_0 \defeq T$. Once $T_n$ is defined, let $\delta_n$ be any group element such that $T_n \cap (T_n\delta_n) = \0$ \ep{such $\delta_n$ exists since $\G$ is infinite} and set $R_n \defeq T_n \sqcup (T_n\delta_n)$. Let $S_n$ be the set $S$ produced by Corollary~\ref{corl:less_technical} applied with $R = R_n$ and $\gamma = \gamma_n$. After replacing $S_n$ with a superset if necessary, we arrange so that $S_n$ is symmetric, contains $\mathbf{1}$, and includes $R_nR_n^{-1}$. Finally, we let $T_{n+1} \defeq S_nT_n$.
        This construction is done so that the following claim holds: %The following claim explains why the sets $T_n$, $R_n$, and $S_n$ are defined in this manner.
	
	\begin{smallclaim}\label{subclaim:split}
		Let $\G \acts X$ be an $S_n$-free continuous action of $\G$ on a zero-dimensional Polish space and let $V \subseteq X$ be a $T_n$-syndetic clopen set. Then there is a partition $V = C \sqcup U$ into two clopen subsets, where $C$ is $R_n$-syndetic, while $U$ is both $S_n$-separated and $T_{n+1}$-syndetic.
	\end{smallclaim}
	\begin{claimproof}[Proof of Claim~\ref{subclaim:split}]
		By Lemma~\ref{lemma:max}, we can let $U$ be a clopen maximal $S_n$-separated subset of $V$ and define $C \defeq V \setminus U$. The maximality of $U$ means that $V \subseteq S_n \cdot U$, and since $V$ is $T_n$-syndetic and $T_{n+1} = S_n T_n$, this implies that $U$ is $T_{n+1}$-syndetic, as claimed.
		To see that $C$ is $R_n$-syndetic, take any $x \in X$. We need to argue that $R_n \cdot x$ contains a point in $C$. Recall that $R_n = T_n \sqcup (T_n \delta_n)$. Since $V$ is $T_n$-syndetic, the sets $T_n \cdot x$ and $T_n \delta_n \cdot x$ each contain a point in $V$. Since the sets $T_n$ and $T_n \delta_n$ are disjoint and the action $\G \acts X$ is $R_nR_n^{-1}$-free, we have $|(R_n \cdot x) \cap V| \geq 2$. On the other hand, $|(R_n \cdot x) \cap U| \leq 1$ since $U$ is $R_nR_n^{-1}$-separated. Therefore, $|(R_n \cdot x) \cap C| \geq 1$, as desired.
	\end{claimproof}
	
	For each $n \in \N$, let $\mathcal{Z}_n \subseteq k^\G$ be the set of all mappings $z \colon \G \to k$ such that
	\[
		\exists \sigma \in S_n \ : \quad z(\sigma) \,\neq\, z(\sigma\gamma_n).
	\]
	The set $\mathcal{Z}_n$ is clopen, and if $z \in \mathcal{Z}_n$, then $\gamma_n \cdot z \neq z$. Hence, the set
	\[
	    \mathcal{Z}^\ast \,\defeq\, \bigcap_{n = 0}^\infty \bigcap_{\delta \in \G} (\delta \cdot \mathcal{Z}_n)
	\]
	is a free subshift. \ep{We have not yet shown that $\mathcal{Z}^\ast$ is nonempty.} For each $N \in \N$, we also let
	\[
		\mathcal{Z}_{< N} \,\defeq\, \bigcap_{n = 0}^{N-1}\bigcap_{\delta \in \G} (\delta \cdot \mathcal{Z}_n).
	\]
	Then $\mathcal{Z}_{< N}$ is a subshift and $\mathcal{Z}^\ast=\dcap_{N=0}^\infty \mathcal{Z}_{< N}$ \ep{where the intersection is decreasing}. %Notice that $Z_{\leq N}$ need not be free; in particular, $\rho$ may not be defined on all of $Z_{\leq N}$.
	We will argue that the conclusion of Lemma~\ref{lemma:extension} holds with $\mathcal{Z} = \mathcal{Z}_{< N}$ for any large enough $N$.
	
	\begin{smallclaim}\label{subclaim:approx}
		For every large enough $N \in \N$, there exists a continuous map $\tilde{\rho} \colon \mathcal{Z}_{< N} \to m$ such that for each $z \in \mathcal{Z}_{< N}$, there is $z^\ast \in \mathcal{Z}^\ast$ with the following properties:
		\begin{enumerate}[label=\ep{\itshape\alph*}]
			\item\label{item:a} for all $\delta \in D$, $z(\delta) = z^\ast(\delta)$; and
			\item\label{item:b} for all $\delta \in D$, $\tilde{\rho}(\delta \cdot z) = \rho(\delta \cdot z^\ast)$.
		\end{enumerate}
	\end{smallclaim}
	\begin{claimproof}[Proof of Claim~\ref{subclaim:approx}]
	    Since $\mathcal{Z}^\ast$ is compact and $\rho$ is continuous, there is a finite set $W \subset \G$ such that for every $z^\ast \in \mathcal{Z}^\ast$, the value $\rho(z^\ast)$ is determined by the restriction of $z^\ast$ to $W$. Observe that for all large enough $N \in \N$ and for each $z \in \mathcal{Z}_{< N}$,
		\begin{equation}\label{eq:approx}
		    \exists z^\ast \in \mathcal{Z}^\ast \ \forall \delta \in D \cup W \cup WD\ : \quad z(\delta) \,=\, z^\ast(\delta).
		\end{equation}
		Indeed, let $\Omega$ be the set of all $z \in k^\G$ for which \eqref{eq:approx} fails. Whether or not $z \in \Omega$ only depends on the restriction of $z$ to $D \cup W \cup WD$, so $\Omega$ is a clopen subset of $k^\G$. Furthermore, $\Omega \cap \mathcal{Z}^\ast = \0$ by definition. Since $k^\G$ is compact and $\mathcal{Z}^\ast=\dcap_{N=0}^\infty \mathcal{Z}_{< N}$, we conclude that $\Omega \cap \mathcal{Z}_{< N} = \0$ for all large enough $N \in \N$, as desired. Now, for any large enough $N$, we can define a function $\tilde{\rho} \colon \mathcal{Z}_{< N} \to k$ by
		\begin{align*}
			\tilde{\rho}(z) = i \quad \vcentcolon&\Longleftrightarrow \quad \exists z^\ast \in \mathcal{Z}^\ast\ : \quad \rho(z^\ast) = i \text{ and } \rest{z^\ast}{W} \,=\, \rest{z}{W}\\
			&\Longleftrightarrow \quad \forall z^\ast \in \mathcal{Z}^\ast\ : \quad \rest{z^\ast}{W} \,=\, \rest{z}{W} \ \Longrightarrow \ \rho(z^\ast) = i.
		\end{align*}
		The two definitions given above are equivalent since for each $z^\ast \in \mathcal{Z}^\ast$, the value $\rho(z^\ast)$ is determined by the restriction of $z^\ast$ to $W$. By construction, $\tilde{\rho}(z)$ is determined by $\rest{z}{W}$, so $\tilde{\rho}$ is continuous. Finally, take any $z \in \mathcal{Z}_{< N}$. By \eqref{eq:approx}, there is $z^\ast \in \mathcal{Z}^\ast$ such that for all $\delta \in D \cup W \cup WD$, $z(\delta) = z^\ast(\delta)$. It is clear that this choice of $z^\ast$ fulfills conditions \ref{item:a} and \ref{item:b}.
	\end{claimproof}
	
	Let $\mathcal{Z} \defeq \mathcal{Z}_{< N}$ for any large enough $N$ and let $\tilde{\rho} \colon \mathcal{Z} \to m$ be given by Claim~\ref{subclaim:approx}. Also, let
	\[
	    S \,\defeq\, \bigcup_{n = 0}^{N-1} S_n.
	\]
	Claim~\ref{subclaim:approx} yields properties \ref{item:Ia} and \ref{item:Ib}, so it remains to verify \ref{item:II}. To this end, let $\G \acts X$ be an $S$-free continuous action of $\G$ on a zero-dimensional Polish space $X$ and let $C_0$, $f_0$ be as in~\ref{item:II}. Since the set $X \setminus C_0$ is $T$-syndetic, we may iteratively apply Claim~\ref{subclaim:split} in order to obtain a sequence of clopen subsets $U_0$, $C_1$, $U_1$, $C_2$, $U_2$, \ldots{} of $X$ such that $U_0 \defeq X \setminus C_0$ and for all $n \in \N$,
	\begin{itemize}
		\item $U_n = C_{n+1} \sqcup U_{n+1}$; and
		\item the set $C_{n+1}$ is $R_n$-syndetic, while $U_{n+1}$ is $S_n$-separated and $T_{n+1}$-syndetic.
	\end{itemize}
	We then use Corollary~\ref{corl:less_technical} repeatedly to obtain an increasing sequence $f_0 \subseteq f_1 \subseteq f_2 \subseteq \ldots$ such that for each $n \in \N$, $f_{n+1} \colon C_0 \sqcup C_1 \sqcup \ldots \sqcup C_{n+1} \to k$ is a continuous function satisfying
	\begin{equation}\label{eq:not_sim}
	    \forall x \in X\ \exists \sigma \in S_n\ : \quad \set{\sigma \cdot x,\, \sigma \gamma_n \cdot x} \,\subseteq\, \dom(f_{n+1}) \text{ and } f_{n+1}(\sigma \cdot x) \,\neq\, f_{n+1}(\sigma\gamma_n \cdot x).
	\end{equation}
	Let $f \colon X \to k$ be an arbitrary continuous extension of $f_{N}$ \ep{for instance, we may set $f(x) \defeq 0$ for all $x \not \in \dom(f_{N})$}. We claim that $f$ is as desired, i.e., that $\pi_f(X) \subseteq \mathcal{Z}$. Indeed, since $\pi_f$ is $\G$-equivariant, it suffices to argue that $\pi_f(x) \in \mathcal{Z}_n$ for all $x \in X$ and $n < N$, i.e., that for all $x \in X$ and $n < N$,
	\[
		\exists \sigma \in S_n\ : \quad  \pi_f(x)(\sigma) \,\neq\, \pi_f(x)(\sigma \gamma_n).
	\]
	Using the definition of $\pi_f$, we can rewrite the latter statement as
	\[
	    \exists \sigma \in S_n\ : \quad f(\sigma \cdot x) \,\neq\, f(\sigma \gamma_n \cdot x),
	\]
	which holds by \eqref{eq:not_sim} since $f$ is an extension of $f_{n+1}$.
	\end{scproof}
	
	\subsubsection{The proof}
	
	In the remainder of this subsection, we prove Lemma~\ref{lemma:cont+Stab}. For convenience, we restate it here: %Let us now state Lemma~\ref{lemma:cont+Stab} again for the reader's convenience:
	
	\begin{lemmacopy}{lemma:cont+Stab}
	    Let $\Sh$ be an SFT and let $\pi \colon \Free(k^\G) \to \Sh$ be a continuous $\G$-equivariant map for some $k \geq 2$. Then there exist a finite symmetric set $F \subset \G$ with $\mathbf{1} \not \in F$, an integer $\ell \geq |F| + 1$, and a continuous $\G$-equivariant map $\tilde{\pi} \colon \Col(F,\ell) \to \Sh$ such that $\Stab(\tilde{\pi}) = \Stab(\pi)$.
	\end{lemmacopy}
	\begin{scproof}\stepcounter{ForClaims} \renewcommand{\theForClaims}{\ref{lemma:cont+Stab}}
	%Let $k$, $\mathcal{S}$, and $\pi$ be as in the statement of Lemma~\ref{lemma:cont+Stab}.
	We may assume that $\pi$ is not constant, since otherwise we can just take $\tilde{\pi}$ to be constant as well. Say $\Sh \subseteq m^\G$ for some $m \geq 1$ and let $\rho \colon \Free(k^\G) \to m$ be the function such that $\pi = \pi_\rho$, i.e.,%. Explicitly, we have
	\[
	    \rho(x) \,\defeq\, \pi(x)(\mathbf{1}) \quad \text{for all } x \in \Free(k^\G).
	\]
	Since $\Sh$ is of finite type, there exist a finite window $W \subset \G$ and a set $\Phi \subseteq m^W$ such that
	\begin{equation}\label{eq:SFT}
	    \Sh \,=\, \set{x \in m^\G \,:\, \rest{(\gamma \cdot x)}{W} \in \Phi \text{ for all } \gamma \in \G}.
	\end{equation}
	Without loss of generality, we may assume that $W$ is symmetric and contains $\mathbf{1}$.
	
	%\begin{remk}
    At this point, it is instructive to notice that Lemma~\ref{lemma:extension} easily yields most of Lemma~\ref{lemma:cont+Stab}:
	
	\begin{smallclaim}\label{claim:noStab}
	    There exist a finite symmetric set $F \subset \G$ with $\mathbf{1} \not \in F$, an integer $\ell \geq |F| + 1$, and a continuous $\G$-equivariant map $\tilde{\pi} \colon \Col(F,\ell) \to \Sh$.
	\end{smallclaim}
	\begin{claimproof}[Proof of Claim~\ref{claim:noStab}]
	    Apply Lemma~\ref{lemma:extension} to $k$, $m$, and $\rho$ and with $D = W$ and $T = \set{\mathbf{1}}$. This yields a subshift $\mathcal{Z} \subseteq k^\G$, a continuous map $\tilde{\rho} \colon \mathcal{Z} \to m$, and a finite set $S \subset \G$ satisfying \ref{item:I} and \ref{item:II}. In particular, since $D = W$, statements \eqref{eq:SFT} and \ref{item:Ib} imply that $\pi_{\tilde{\rho}}(\mathcal{Z}) \subseteq \Sh$. Without loss of generality, we may assume that $S$ is symmetric and contains $\mathbf{1}$. Take \[F \,\defeq\, S \setminus \set{\mathbf{1}} \quad \text{and} \quad \ell \,\defeq\, |F| + 1.\] The shift action $\G \acts \Col(F, \ell)$ is $S$-free, so we may apply part \ref{item:II} of Lemma~\ref{lemma:extension} to $X = \Col(F, \ell)$ with $C_0 = \0$ and $f_0 = \0$. This yields a continuous map $f \colon \Col(F, \ell) \to k$ with $\pi_f(\Col(F, \ell)) \subseteq \mathcal{Z}$. The composition $\tilde{\pi} \defeq \pi_{\tilde{\rho}} \circ \pi_f$ is a continuous $\G$-equivariant map from $\Col(F, \ell)$ to $\Sh$, as desired.
	\end{claimproof}
	
	The only issue with Claim~\ref{claim:noStab} is that it gives no control over $\Stab(\tilde{\pi})$. To ensure that $\Stab(\tilde{\pi}) = \Stab(\pi)$, we shall invoke the same construction, but starting with a carefully chosen function $f_0$. 
	%\end{remk}
	%After we prove Claim~\ref{claim:extension}, we shall apply its part \ref{item:II} with $Y = \Col(F, \ell)$ for appropriately chosen $F$ and $\ell$. This will give us a continuous map $f \colon \Col(F, \ell) \to k$ such that $\pi_f(\Col(F, \ell)) \subseteq Z$. Assuming that $W \subseteq D$, \eqref{eq:SFT} and property \ref{item:Ib} imply that $\pi_{\tilde{\rho}}(Z) \subseteq X$, so the composition $\tilde{\pi} \defeq \pi_{\tilde{\rho}} \circ \pi_f$ will be a continuous $\G$-equivariant map from $\Col(F, \ell)$ to $X$. Finally, the sets $D$, $T$ and the function $f_0$ will be carefully constructed to ensure that $\Stab(\tilde{\pi}) = \Stab(\pi)$.
	
	Our argument proceeds in several stages.
	
	\medskip
	        
	\noindent\textbf{Stage 1:} \textsl{Constructing sets $D$ and $T$}. 
	To begin with, we define the sets $D$ and $T$ to which Lemma~\ref{lemma:extension} will be applied. %The set $D$ is built in two steps. First,
	Since $\rho$ is not constant, there are two points $a_0$, $a_1 \in \Free(k^\G)$ with \[\rho(a_0) \,\neq\, \rho(a_1).\] Using the continuity of $\rho$, we obtain a finite set $A \subset \G$ such that for all $x \in \Free(k^\G)$, $\rho(x) = \rho(a_0)$ or $\rho(x) = \rho(a_1)$ whenever $x$ agrees with $a_0$ or $a_1$ respectively on $A$. Without loss of generality, we may assume that the set $A$ is symmetric and contains $\mathbf{1}$.
	
	Next, %let $T_0$ be an arbitrary finite symmetric subset of $\G$ with $\mathbf{1} \in T_0$ such that $|T_0| > |A|$. Set \[r \,\defeq\, |AT_0^2A \setminus \Stab(\pi)|,\] and let $\gamma_1$, \ldots, $\gamma_r$ be an enumeration of the elements of $AT_0^2A \setminus \Stab(\pi)$. For each $1 \leq i \leq r$, since $\gamma_i \not \in \Stab(\pi)$, we can pick $b_i \in \Free(k^\G)$ and $\sigma_i \in \G$ so that
	%\[
	%    \rho(\sigma_i \cdot b_i) \,\neq\, \rho(\sigma_i \gamma_i \cdot b_i).
	%\]
	%Since $\rho$ is continuous, there is a finite set $B_i \subset \G$ such that for all $x \in \Free(k^\G)$,
	%\[
	%    \rest{x}{B_i} \,=\, \rest{b_i}{B_i} \quad \Longrightarrow \quad \rho(\sigma_i \cdot x) \,=\, \rho(\sigma_i \cdot b_i) \text{ and } \rho(\sigma_i \gamma_i \cdot x) \,=\, \rho(\sigma_i \gamma_i \cdot b_i).
	%\]
	%In particular, every point $x \in \Free(k^\G)$ that agrees with $b_i$ on $B_i$ satisfies $\rho(\sigma_i \cdot x) \neq \rho(\sigma_i \gamma_i \cdot x)$. 
	given any group element $\gamma \not\in \Stab(\pi)$, we can pick some $b_\gamma \in \Free(k^\G)$ and $\sigma_\gamma \in \G$ so that
	\[
	    \rho(\sigma_\gamma \cdot b_\gamma) \,\neq\, \rho(\sigma_\gamma \gamma \cdot b_\gamma).
	\]
	Since $\rho$ is continuous, there is a finite set $B_\gamma \subset \G$ such that for all $x \in \Free(k^\G)$,
	\[
	    \rest{x}{B_\gamma} \,=\, \rest{b_\gamma}{B_\gamma} \quad \Longrightarrow \quad \rho(\sigma_\gamma \cdot x) \,=\, \rho(\sigma_\gamma \cdot b_\gamma) \text{ and } \rho(\sigma_\gamma \gamma \cdot x) \,=\, \rho(\sigma_\gamma \gamma \cdot b_\gamma).
	\]
	In particular, every point $x \in \Free(k^\G)$ that agrees with $b_\gamma$ on $B_\gamma$ satisfies \[\rho(\sigma_\gamma \cdot x) \,\neq\, \rho(\sigma_\gamma \gamma \cdot x).\]
	By replacing each $B_\gamma$ with a superset if necessary, we may assume that $B_\gamma$ is symmetric and
	\[
	    \set{\mathbf{1}, \, \sigma_\gamma, \, \sigma_\gamma \gamma} \, \subseteq\, B_\gamma.
	\]
	
	Now we let $T_0$ be an arbitrary finite symmetric subset of $\G$ with $\mathbf{1} \in T_0$ and $|T_0| > |A|$. Define
	\[
	    B \,\defeq\, \bigcup \big\{B_\gamma \,:\, \gamma \in AT_0^2A\setminus \Stab(\pi)\big\}.
	\]
	Let $D \defeq W \cup A  \cup B \cup \Stab(\pi)$. \ep{Note that $\Stab(\pi)$ is finite by Proposition~\ref{prop:Stab}, so $D$ is finite as well.} %Note that $D$ is symmetric and contains $\mathbf{1}$.
	Then we let $T_1$ be an arbitrary finite symmetric subset of $\G$ such that $\mathbf{1} \in T_1$ and $|T_1| > |T_0||B|$, and define $T \defeq T_0T_1$.

	\medskip
	        
	\noindent\textbf{Stage 2:} \textsl{Defining $F$ and $\ell$}.  Applying Lemma~\ref{lemma:extension} to the sets $D$ and $T$ from Stage $1$ yields
	\begin{itemize}
	        \item a subshift $\mathcal{Z} \subseteq k^\G$,
	        \item a continuous map $\tilde{\rho} \colon \mathcal{Z} \to m$, and
	        \item a finite set $S \subset \G$
	\end{itemize}
	%\[
	%    \text{a subshift $Z \subseteq k^\G$}, \quad \text{a continuous map $\tilde{\rho} \colon Z \to m$}, \quad \text{and} \quad \text{a finite set $S \subset \G$}
	%\]
	%\begin{itemize}
	%        \item a subshift $Z \subseteq k^\G$,
	%        \item a continuous map $\tilde{\rho} \colon Z \to m$, and
	%        \item a finite set $S \subset \G$
	%\end{itemize}
	satisfying conditions \ref{item:I} and \ref{item:II}. We may assume that $S$ is symmetric and $\mathbf{1} \in S$. Let
	\[
	    F \,\defeq\, (T_0 \cup T_1 \cup D \cup S)^{100}\setminus \set{\mathbf{1}} \qquad \text{and} \qquad \ell \defeq |F| + 4.
	\]
	Note that the shift action $\G \acts \Col(F, \ell)$ is $F$-free. In particular, it is $S$-free, so Lemma~\ref{lemma:extension}\ref{item:II} can be applied to it. We will show that the conclusion of Lemma~\ref{lemma:cont+Stab} holds with this choice of $F$ and $\ell$.
	
	\medskip
	        
	\noindent\textbf{Stage 3:} \textsl{Constructing $C_0$ and $f_0$}. For $i \in \set{0,1}$, let $J_i$ be the clopen subset of $\Col(F, \ell)$ given by
	\[
	    J_i \,\defeq\, \set{x \in \Col(F, \ell) \,:\, x(\mathbf{1}) = i \text{ and } x(\delta) \geq 2 \text{ for all } \delta \in AT_0^2A \setminus \set{\mathbf{1}}}.
	\]
	Note that $J_0 \cap J_1 = \0$. Moreover, the union $J_0 \sqcup J_1$ is $AT_0^2A$-separated. In particular, the sets $A \cdot J_0$ and $A \cdot J_1$ are disjoint, so we can set \[\mathfrak{A} \,\defeq\, (A \cdot J_0) \sqcup (A \cdot J_1).\]
	Since the shift action $\G \acts \Col(F, \ell)$ is $A^2$-free, for each point $y \in \mathfrak{A}$, there is precisely one choice of $x \in J_0 \sqcup J_1$ and $\alpha \in A$ such that $y = \alpha \cdot x$. Thus, we can define a continuous function $g \colon \mathfrak{A} \to k$ by %setting, for each $\delta \in A$ and $x \in J_0 \sqcup J_1$,
	\[
	    g(\alpha \cdot x) \,\defeq\,
	        a_i(\alpha) \quad \text{for all } x \in J_i \text{ and } \alpha \in A.
	\]
	In other words, $g$ is obtained by copying $\rest{a_i}{A}$ to $A \cdot x$ for all $x \in J_i$. %\ep{Here we are using that the action $\G \acts \Col(F, \ell)$ is $F$-free and $F \supseteq A^2\setminus \set{\mathbf{1}}$.}
	
	Using Lemma~\ref{lemma:max}, we let $K$ be a clopen maximal $BTT^{-1}B$-separated subset of $\Col(F, \ell) \setminus (B \cdot \mathfrak{A})$. \ep{Here we are using that the shift action $\G \acts \Col(F, \ell)$ is $BTT^{-1}B$-free.} 
	
	\begin{smallclaim}\label{claim:K}
	    The set $K$ is infinite.
	\end{smallclaim}
	\begin{claimproof}[Proof of Claim~\ref{claim:K}]
	    It is easy to construct infinitely many distinct proper colorings of the Cayley graph $G(\G, F)$ using the colors $2$, $3$, \ldots, $\ell - 1$. For instance, one can assign the color $2$ to an arbitrary vertex of $G(\G, F)$ and then use only the colors $3$, \ldots, $\ell - 1$ on the remaining vertices \ep{this is possible since the list $3$, \ldots, $\ell - 1$ includes $\ell - 3 = |F| + 1$ colors}. If $x \in \Col(F, \ell)$ is a coloring that only uses the colors $2$, $3$, \ldots, $\ell - 1$, then the orbit of $x$ under the shift action $\G \acts \Col(F, \ell)$ is disjoint from $\mathfrak{A}$, and, in particular, $x \not \in B \cdot \mathfrak{A}$. Thus, the set $\Col(F, \ell) \setminus (B \cdot \mathfrak{A})$ is infinite. Since $\Col(F, \ell) \setminus (B \cdot \mathfrak{A}) \subseteq BTT^{-1}B \cdot K$, this implies that $K$ must be infinite as well.
	\end{claimproof}
	
	Since the set $K$ is infinite, we may partition it as \[K \,=\, \bigsqcup_\gamma K_\gamma,\] where the union is over all $\gamma \in AT_0^2A\setminus \Stab(\pi)$ and each $K_\gamma$ is nonempty and clopen. Set $\mathfrak{B} \defeq B \cdot K$. Then for each point $y \in \mathfrak{B}$, there is precisely one choice of $x \in K$ and $\beta \in B$ such that $y = \beta \cdot x$. Therefore,  %Note that  Let \[\mathscr{B}_\gamma \,\defeq\, B_\gamma \cdot K_\gamma \quad \text{and} \quad \mathscr{B} \,\defeq\, \bigsqcup_\gamma \mathscr{B}_\gamma\] \ep{the union is disjoint since $K$ is $D^2$-separated} and 
	we can define a continuous function $h \colon \mathfrak{B} \to k$ by
	\begin{equation}\label{eq:b}
		h(\beta \cdot x) \,\defeq\, b_\gamma(\beta) \quad \text{for all } x \in K_\gamma \text{ and } \beta \in B.
	\end{equation}
	In other words, we copy $\rest{b_\gamma}{B}$ onto $B \cdot x$ for each $x \in K_\gamma$.
	%Property \eqref{eq:b} will be eventually used to show that $\Pat_F(X, g) \supseteq \Pat_F(\Free(2^\G), f)$.
	
	Since $K \subseteq \Col(F, \ell) \setminus (B \cdot \mathfrak{A})$, the sets $\mathfrak{A}$ and $\mathfrak{B}$ are disjoint. Thus, we can take the disjoint union $C_0 \defeq \mathfrak{A} \sqcup \mathfrak{B}$ and let $f_0 \colon C_0 \to k$ be the function that is equal to $g$ on $\mathfrak{A}$ and $h$ on $\mathfrak{B}$.
	
	\medskip
	        
	\noindent\textbf{Stage 4:} \textsl{Using Lemma~\ref{lemma:extension}}. To apply part \ref{item:II} of Lemma~\ref{lemma:extension} to the function $f_0$ constructed in Stage 3, we need to verify that the set $\Col(F, \ell) \setminus C_0$ is $T$-syndetic.
	
	\begin{smallclaim}\label{claim:T}
	    The set $\Col(F, \ell) \setminus C_0$ is $T$-syndetic.
	\end{smallclaim}
	\begin{claimproof}[Proof of Claim~\ref{claim:T}]
	    We first show that the set $\Col(F, \ell) \setminus \mathfrak{A}$ is $T_0$-syndetic. To this end, take any $x \in \Col(F, \ell)$. We have to show that $(T_0 \cdot x) \setminus \mathfrak{A} \neq \0$. Note that $|T_0 \cdot x| = |T_0|$, since the shift action $\G \acts \Col(F, \ell)$ is $T_0$-free. By construction, the set $J_0 \sqcup J_1$ is $AT_0^2A$-separated, so there is at most one point $y \in J_0 \sqcup J_1$ with $(T_0 \cdot x) \cap (A \cdot y) \neq \0$. Hence, $|(T_0 \cdot x) \setminus \mathfrak{A}| \geq |T_0| - |A| > 0$, as desired.
	    
	    Now we need to show that for every $x \in \Col(F,\ell)$, $(T \cdot x) \setminus C_0 \neq \0$. Recall that $T = T_0T_1$. Since the set $\Col(F, \ell) \setminus \mathfrak{A}$ is $T_0$-syndetic, $(T_0 \tau \cdot x) \setminus \mathfrak{A} \neq \0$ for all $\tau \in T_1$. Therefore,
	    \begin{equation}\label{eq:q}
	        |(T \cdot x) \setminus \mathfrak{A}| \,\geq\, \frac{|T_1|}{|T_0|} \,>\, |B|.
	    \end{equation}
	    The set $K$ is $BTT^{-1}B$-separated, so there is at most one $y \in K$ with $(T \cdot x) \cap (B \cdot y) \neq \0$. Hence,
	    \begin{equation}\label{eq:w}
	        |(T \cdot x) \cap \mathfrak{B}| \,\leq\, |B|.
	    \end{equation}
	    From \eqref{eq:q} and \eqref{eq:w}, it follows that $|(T \cdot x) \setminus C_0| > 0$, and we are done.
	\end{claimproof}
	
	Thanks to Claim~\ref{claim:T} and since the action $\G \acts \Col(F, \ell)$ is $S$-free, we may apply Lemma~\ref{lemma:extension}\ref{item:II} and obtain a continuous function $f \colon \Col(F, \ell) \to k$ that extends $f_0$ and satisfies
	\[
	    \pi_f(\Col(F, \ell)) \,\subseteq\, \mathcal{Z}.
	\]
	We claim that the composition $\tilde{\pi} \defeq \pi_{\tilde{\rho}} \circ \pi_f$ satisfies the conclusion of Lemma~\ref{lemma:cont+Stab}.
	
	\medskip
	        
	\noindent\textbf{Stage 5:} \textsl{Finishing the proof}. Since $D \supseteq W$, statements \eqref{eq:SFT} and \ref{item:Ib} imply that $\pi_{\tilde{\rho}}(\mathcal{Z}) \subseteq \Sh$. Thus, $\tilde{\pi}$ is a continuous $\G$-equivariant map from $\Col(F, \ell)$ to $\Sh$. It remains to verify that $\Stab(\tilde{\pi}) = \Stab(\pi)$.
	
	We start with the easier inclusion $\Stab(\tilde{\pi}) \supseteq \Stab(\pi)$:
	
	\begin{smallclaim}\label{claim:more}
	    $\Stab(\tilde{\pi}) \supseteq \Stab(\pi)$.
	\end{smallclaim}
	\begin{claimproof}[Proof of Claim~\ref{claim:more}]
	    It suffices to show that $\Stab(\pi_{\tilde{\rho}}) \supseteq \Stab(\pi)$. To this end, we first claim that for all $z \in \mathcal{Z}$ and $\gamma \in \Stab(\pi)$, $\tilde{\rho}(z) = \tilde{\rho}(\gamma \cdot z)$. Indeed, since $D \supseteq \Stab(\pi)$, by part \ref{item:Ib} of Lemma~\ref{lemma:extension}, there is a point $z^\ast \in \Free(k^\G)$ such that $\tilde{\rho}(z) = \rho(z^\ast)$ and $\tilde{\rho}(\gamma \cdot z) = \rho(\gamma \cdot z^\ast)$. But $\gamma \in \Stab(\pi(z^\ast))$, so in particular $\rho(z^\ast) = \rho(\gamma \cdot z^\ast)$, as claimed. Now take any $z \in \mathcal{Z}$ and $\gamma \in \Stab(\pi)$. Since $\Stab(\pi)$ is a normal subgroup of $\G$, for each $\sigma \in \G$, we have
	    \[
	        \tilde{\rho}(\sigma \cdot z) \,=\, \tilde{\rho}(\sigma \gamma \sigma^{-1} \cdot (\sigma \cdot z)) \,=\, \tilde{\rho}(\sigma \gamma \cdot z).
	    \]
	    This means that $\pi_{\tilde{\rho}}(z) = \pi_{\tilde{\rho}}(\gamma \cdot z)$, i.e., $\gamma \in \Stab(\pi_{\tilde{\rho}}(z))$, as desired.
	\end{claimproof}
	
	In the next two claims we reap the fruits of the labor we invested in the construction of $f_0$. 
	
	\begin{smallclaim}\label{claim:A}
	    If $x \in J_i$ for some $i \in \set{0,1}$, then $\tilde{\rho}(\pi_f(x)) = \rho(a_i)$.
	\end{smallclaim}
	\begin{claimproof}[Proof of Claim \ref{claim:A}]
	    Set $z \defeq \pi_f(x)$. Since $D \supseteq A$, Lemma~\ref{lemma:extension}\ref{item:I} yields a point $z^\ast \in \Free(k^\G)$ such that $\tilde{\rho}(z) = \rho(z^\ast)$ and $\rest{z}{A} = \rest{z^\ast}{A}$. As $x \in J_i$ and $f$ extends $g$, we have $\rest{z}{A} = \rest{a_i}{A}$. %by construction.
	    Therefore, $z^\ast$ agrees with $a_i$ on $A$ and hence, by the choice of $A$, $\rho(z^\ast) = \rho(a_i)$, as desired.
	\end{claimproof}
	
	\begin{smallclaim}\label{claim:B}
	    If $x \in K_\gamma$ for some $\gamma \in AT_0^2A\setminus \Stab(\pi)$, then $\tilde{\rho}(\sigma_\gamma \cdot \pi_f(x)) \neq \tilde{\rho}(\sigma_\gamma \gamma \cdot \pi_f(x))$.
	\end{smallclaim}
	\begin{claimproof}[Proof of Claim \ref{claim:B}]
	    Set $z \defeq \pi_f(x)$. Since $D \supseteq B \supseteq B_\gamma \supseteq \set{\sigma_\gamma, \sigma_\gamma \gamma}$, Lemma~\ref{lemma:extension}\ref{item:I} gives a point $z^\ast \in \Free(k^\G)$ such that $\tilde{\rho}(\sigma_\gamma \cdot z) = \rho(\sigma_\gamma \cdot z^\ast)$, $\tilde{\rho}(\sigma_\gamma \gamma \cdot z) = \rho(\sigma_\gamma \gamma \cdot z^\ast)$, and $\rest{z}{B_\gamma} = \rest{z^\ast}{B_\gamma}$. %with the following properties:
	    %\begin{itemize}
	    %    \item $\tilde{\rho}(\sigma_\gamma \cdot z) = \rho(\sigma_\gamma \cdot z^\ast)$,
	    %    \item $\tilde{\rho}(\sigma_\gamma \gamma \cdot z) = \rho(\sigma_\gamma \gamma \cdot z^\ast)$, and
	    %    \item $\rest{z}{B_\gamma} = \rest{z^\ast}{B_\gamma}$.
	    %\end{itemize}
	    Since $x \in K_\gamma$ and $f$ extends $h$, we have $\rest{z}{B_\gamma} = \rest{b_\gamma}{B_\gamma}$. Therefore, $z^\ast$ agrees with $b_\gamma$ on $B_\gamma$ and hence, by the choice of $B_\gamma$, $\rho(\sigma_\gamma \cdot z^\ast) \neq \rho(\sigma_\gamma \gamma \cdot z^\ast)$, as desired.
	\end{claimproof}
	
	Finally, we are ready to establish the inclusion $\Stab(\tilde{\pi}) \subseteq \Stab(\pi)$:
	
	\begin{smallclaim}\label{claim:less}
	    $\Stab(\tilde{\pi}) \subseteq \Stab(\pi)$.
	\end{smallclaim}
	\begin{claimproof}[Proof of Claim~\ref{claim:less}]
	    Take any group element $\gamma \not \in \Stab(\pi)$. We have to find a point $x \in \Col(F, \ell)$ such that $\tilde{\pi}(x) \neq \tilde{\pi}(\gamma \cdot x)$. We consider two cases.
	    
	    \smallskip
	        
	    \noindent \textbf{Case 1:} \textsl{$\gamma \not \in AT_0^2A$}. Since $\gamma \neq \mathbf{1}$, the mapping $\mathbf{1} \mapsto 0$, $\gamma \mapsto 1$ is a proper partial coloring of the Cayley graph $G(\G, F)$. The list $2$, $3$, \ldots, $\ell - 1$ includes $\ell - 2 > |F| + 1$ colors, so this mapping can be extended to a coloring $x \in \Col(F, \ell)$ such that $x(\delta) \geq 2$ for all $\delta \not \in \set{\mathbf{1}, \gamma}$. Since $\gamma \not \in AT_0^2A$, we have $x \in J_0$, while $\gamma \cdot x \in J_1$. By Claim~\ref{claim:A}, \[\tilde{\rho}(\pi_f(x)) \,=\, \rho(a_0) \,\neq\, \rho(a_1) \,=\, \tilde{\rho}(\pi_f(\gamma \cdot x)),\] i.e., $\tilde{\pi}(x)(\mathbf{1}) \neq \tilde{\pi}(\gamma \cdot x)(\mathbf{1})$. Therefore, $\tilde{\pi}(x) \neq \tilde{\pi}(\gamma \cdot x)$, as desired.
	    
	    \smallskip
	        
	    \noindent \textbf{Case 2:} \textsl{$\gamma \in AT_0^2A \setminus \Stab(\pi)$}. Take any $x \in K_\gamma$. By Claim~\ref{claim:B}, we have
	    \[
	        \tilde{\rho}(\sigma_\gamma \cdot \pi_f(x)) \,\neq\, \tilde{\rho}(\sigma_\gamma \gamma \cdot \pi_f(x)),
	    \]
	    i.e., $\tilde{\pi}(x)(\sigma_\gamma) \neq \tilde{\pi}(\gamma \cdot x)(\sigma_\gamma)$. Therefore, $\tilde{\pi}(x) \neq \tilde{\pi}(\gamma \cdot x)$, and we are done.
	\end{claimproof}
	
	Together Claims~\ref{claim:more} and \ref{claim:less} imply that \[\Stab(\tilde{\pi}) \,=\, \Stab(\pi),\] which completes the proof of Lemma~\ref{lemma:cont+Stab}.
	\end{scproof}
	
	%\section{Proof of Theorem~\ref{theo:main}}\label{sec:proof}
	
	%\subsection{First stage}\label{subsec:first_stage}
	
	%\subsection{Second stage}\label{subsec:second_stage}
	
	%\subsection{Third stage}\label{subsec:third_stage}
	
	%\subsection{Set-up}
	
	\printbibliography
    
\end{document}